\documentclass[final,leqno,onefignum,onetabnum]{siamltexmm}

\usepackage[letterpaper, margin=1in]{geometry}

\usepackage{amsmath}
\usepackage{amsfonts}
\usepackage{amssymb}

\usepackage{circuitikz}
\usepackage{graphicx,xcolor}

\title{Numerical Diagnostics for Systems of Differential Algebraic Equations} 

\author{Matthew O. Williams\thanks{United Technologies Research Center, 411 Silver Lane, East Hartford, CT 06108 USA.
 (\email{william1@utrc.utc.com})} \and Teems E. Lovett\thanks{United Technologies Research Center, 411 Silver Lane, East Hartford, CT 06108 USA. (\email{LovettTE@utrc.utc.com})}}

\begin{document}
\maketitle
\newcommand{\slugmaster}{%
\slugger{siads}{xxxx}{xx}{x}{x--x}}%slugger should be set to juq, siads, sifin, or siims

\begin{abstract}
In many commercial and academic settings, numerical solvers fail to achieve  their theoretical performance levels due to issues in the system definition, parameterization, and even implementation.
We propose a pair of methods for detecting and localizing these convergence rate issues in applications that can be treated as homotopy problems including  numerical continuation and the evolution of differential algebraic equations. 
Both approaches are rooted in dynamical systems theory, in particular, the numerical techniques used to perform bifurcation studies on ``black-box'' systems, and can be applied across a range of numerical solvers and systems without significant modification.
In a general problem, these methods will detect certain classes of  {\em convergence anomalies}, and indicate the states that are affected by their existence. 
However in applications like circuit simulation, certain classes of anomalies can be localized to an individual component through the use of a small number of additional computational experiments.
To demonstrate the efficacy of the approaches for detection and localization, we apply them to a pair of nonlinear circuits: a diode bridge and a simplified model of a power channel.
\end{abstract}



\pagestyle{myheadings}
\thispagestyle{plain}
\markboth{M.O. Williams \& T.E. Lovett}{Numerical Diagnostics for DAEs}

\section{Introduction}
\label{sec:intro}

Whether as part of the dynamical simulation of a time-evolving system~\cite{shampine1979user}, the optimization of a function in the presence of constraints~\cite{potra2000interior,biegler2009large}, the tracking of steady state operating points as a parameter changes~\cite{doedel1981auto,salinger2011loca,dhooge2003matcont}, or numerous other applications, solving nonlinear systems of equations is an integral part of many numerical routines.
Although the theory behind these methods is well established~\cite{nocedal2006numerical,kelley2003solving}, solvers frequently fail to achieve the convergence rates they are theoretically capable of achieving. 
The root causes of these failures vary immensely and includes issues as localized as a ``typo'' in the system implementation to ones as diffuse as a less than ideal parameterization of the solver and/or system.
Ultimately, the goal of {\em numerical diagnostics} is to identify these issues when they appear, and to suggest corrective actions for the developers of numerical solvers or the systems they are applied to.
This is particularly important in commercial settings where the system state can be large and the implementation hidden in a ``black box'' that prevents detailed analysis of its inner workings.

In response to this need, a number of different approaches have been proposed. 
Most commonly, techniques like Newton's method are augmented with line-searches or trust regions to make the method more robust to poor initial guesses or other issues~\cite{nocedal2006numerical,kelley2003solving}.
The objective of these approaches is to enable a single solver to perform well on a larger class of problems, hence reducing the probability that numerical issues are encountered.
In the context of optimization, recent work by Drori and Teboulle~\cite{drori2014performance} and Lessard, Recht, and Packard~\cite{lessard2014analysis} have resulted in methods with provable upper bounds on their convergence rates.
Because solving nonlinear systems of equations can be recast as an optimization problem, these results can be employed more generally, and their convergence rate bounds provide confidence in the validity and accuracy of the resulting solution.
However, they are only applicable for first order optimization methods, and many commercial or scientific codes~\cite{keiter2013xyce,biegler2009large,dhooge2003matcont,salinger2011loca,doedel1981auto} use other approaches.

In this manuscript, we present a pair of methods that can be applied to a broad class of nonlinear solvers, but assume the underlying problem is a homotopy problem.
In addition to applications such as {\em gmin homotopy}~\cite{keiter2013xyce,vladimirescu1994spice}, which is used to find DC operating points in electronic circuits, other tasks such as time stepping and numerical continuation~\cite{doedel1981auto,dhooge2003matcont} can also be treated as homotopy problems.
The free parameters in those applications can be thought of as time and arclength respectively, and unlike ``standard'' homotopy problems, we are interested in the intermediate solutions as well as the solution at the initial and terminal times.
In general,  our approaches can be applied when: (i) the problem can be reduced to a sequence of nonlinear root finding problems where (ii) a sequence of good initial guesses are available from previously obtained data.

The theory behind these approaches is based on standard dynamical systems techniques for bifurcation studies.
Treating a nonlinear solver as a dynamical system is not new, and in fact, often one step in the associated convergence proof~\cite{nocedal2006numerical}.
However, the solvers and systems that are implemented in practice often differ from their theoretical analogs.
This mismatch between theory and practice could be due to a deliberate design choice, such as when Newton's method is replaced with Newton-Krylov-GMRES~\cite{kelley2003solving}, or could be unintentional as is the case with  implementation errors in the underlying system.
To capture the pairing of system and solver as they are implemented,  we adapt techniques used to perform bifurcation studies on black-box  systems~\cite{lehoucq1998arpack,kelley2003solving,shroff1993stabilization,kevrekidis2003equation}, which are by now well established even in problems with high-dimensional state spaces.
Therefore, our contribution is not new analysis or numerical methods, but the combination of the analyses commonly performed on idealized solvers with techniques that can be applied to practical systems.
As a result, the methods we proposed enable the detection of {\em convergence anomalies}, which is when the actual convergence rate of the solver deviates significantly from the theoretical baseline, and depending on the level of information available and the specific problem identified, the possibility of isolating these issues to the equations responsible for the deviation between theory and practice.

In the remainder of this manuscript, we will discuss our dynamical systems based approaches to numerical diagnostics and demonstrate them on a pair of illustrative examples.
To provide further information about the class of problems these methods can be applied to, we will go into further detail about our definition of a homotopy problem in Sec.~\ref{sec:homotopy}.
In Sec.~\ref{sec:dynamical-systems} we review how the results of linear stability analysis can help diagnose system issues, and provide our approaches for performing this analysis in practice.
Although two approaches will be presented, the main difference between them is in their implementation: the first approach uses standard techniques for linear stability analysis in large systems, while the second approach is data-driven and meant to be used when commercial software packages, which often do not allow the level of control needed for the first approach, are used. 
In Sec~\ref{sec:circuits}, we demonstrate these methods on examples taken from the simulation of circuits, and show how additional knowledge about these systems can isolate implementation errors down to individual electrical components.
Finally, in Sec.~\ref{sec:conclusions}, we have some brief concluding remarks and future outlook.

\section{Homotopy Problems}
\label{sec:homotopy}

There are a number of issues that could arise while solving nonlinear systems of equations, and one of the most prevalent is when techniques like Newton's method are supplied with a ``bad initial guess.''
Because Newton's method is infamous for having fractal basins on attraction, performing a full nonlinear analysis to determine if and when the methods converges when applied to an arbitrary black box system is well beyond the scope of this manuscript.
To circumvent the issue of initialization, we restrict our focus to applications that can be treated as homotopy problems.

The quintessential homotopy problem is of the form:
\begin{equation}
\label{eq:homotopy:quint}
h(x, \theta) = \theta f(x) + (1-\theta)g(x) = 0,
\end{equation}
where $f$ is the function whose roots are desired, $g$ is an auxiliary function whose roots are known, and $\theta$ is a parameter that will be slowly varied from $0$ to $1$. 
To compute $x$ such that $f(x) = 0$, we solve a sequence of nonlinear systems $h(x, \theta) = 0$ by slowly increasing $\theta$ from zero to one.
At $\theta = 0$, the system of interest is $h(x, 0) = g(x) = 0$, whose solution is assumed to be known  (or relatively simple to compute).
By stepping in $\theta$, the previous solution can be used to initialize the nonlinear solver, and provides the ``good initial guess'' these methods often require.
The procedure terminates at $\theta=1$, where $h(x, 1) = f(x)  = 0$, which is the problem whose solution is desired.

Homotopy problems appear in a number of contexts.
Computationally, numerical continuation packages such as LOCA~\cite{salinger2011loca} implement \eqref{eq:homotopy:quint} directly as one method for solving difficult systems of nonlinear equations.  
Interior point methods, which implement inequality constraints through the use of logarithmic barrier functions whose coefficients are slowly reduced, can also be recast as homotopy problems~\cite{potra2000interior}. 
In circuit simulation, they are frequently used to identify DC operating points in systems with nonlinear components~\cite{keiter2013xyce,vladimirescu1994spice, melville1993artificial,yamamura1999fixed}, and are implemented in various ways including the slow removal of artificially introduced pathways to ground or ``sweeping''  the operating voltage from zero to some desired value.

The ``homotopy problem'' that will appear in our examples is the time-integration of ODEs or DAEs using implicit integrators like the Backward Differentiation Formulas (BDFs)~\cite{shampine1979user}. 
Each time step taken using a BDF with the fixed time step $\Delta t$ requires the solution of the nonlinear system 
\begin{equation}
\label{eq:homotopy:time step}
F(t, x, \dot x) = 0.
\end{equation}
To convert this system to an algebraic system of equations, we approximate $\dot x = \alpha x - \beta$, where $\alpha$ is a coefficient that depends on the order of the BDF and the time step, and $\beta$ is a vector containing a linear combination of the solutions at previous time steps.
If the governing equations are the ODE $\dot x = f(t, x)$ and a first order BDF is used, then $F$ in (\ref{eq:homotopy:time step}) is 
\begin{equation}
\label{eq:homotopy:simple-ode}
F(t, x, \beta) =   \alpha x - \beta - f(t+\Delta t, x) = 0,
\end{equation}
where $\alpha = 1/\Delta t$ and $\beta$ is the solution at the previous time step scaled by $1/\Delta t$. 
As stated previously,  this problem has two features that allows our approach to be of use in this setting.
The first is that ``good initial guesses'' for the nonlinear solver can be obtained using the previous system state or some extrapolation based on that data, which circumvents the need to perform a full nonlinear analysis of the system.
The second is that evolving this system in time requires multiple solves of \eqref{eq:homotopy:time step}, which is a restriction whose importance will be discussed in the next section.

\section{Dynamical Systems Approach to Numerical Diagnostics}
\label{sec:dynamical-systems}

In this section, we demonstrate how approaches taken from dynamical systems can be adapted to diagnose numerical problems that arise when solving homotopy problems such as dynamically simulating nonlinear systems. 
To outline our approach, we will first review how the solution of a nonlinear system of equations can be rewritten as a nonlinear, discrete-time dynamical system, and how linear stability analysis can be used to characterize the effective convergence rate of a nonlinear solver, i.e., the performance of a nonlinear solver {\em as it is implemented}, and identify deviations from the theoretically predicted baseline.
In particular, we are interested in {\em convergence anomalies}, which we define as the convergence rate-limiting behaviors that arise in the nonlinear solving process.
While there are many possible causes of a convergence  anomaly, we propose a directional-derivative based check that can often localize the source of the anomaly to a subset of the residual equations.
Finally, because our interest is in diagnosing numerical issues in homotopy problems, we will present two approaches for tracking convergence anomalies as the homotopy parameter changes: the first requires the ability to run additional computational experiments, while the second is a purely data-driven approach.

\subsection{Nonlinear Solvers as Iterative Maps}

In many texts, Newton's method is presented as a two-step procedure~\cite{press1996numerical}.
Given an initial guess, $x_n$, one first computes a ``correction'' by solving a linear system $\delta x_n = J^{-1}(x_n)F(x_n)$, where $x_n$ is the guess at the $n$-th iterate, $F$ is the residual function whose roots we are interested in, and $J = \frac{\partial F}{\partial x}$ is the Jacobian of $F$. 
With this correction, the guess is now updated $x_{n+1} = x_n - \delta x_n$, and the procedure continues until the residual of $F$ is below some predetermined tolerance.
However,  Newton's method can also be written as a dynamical system with the evolution law:
\begin{equation}
\label{eq:theory:newton}
x_{n+1} = x_n - J^{-1}(x_n)F(x_n),
\end{equation}
where the objective is now to find a fixed point of this map.
Provided that $F$ has an isolated root, continuous second derivatives, and that we have a good enough initial guess, Newton's method as presented in (\ref{eq:theory:newton}) will converge to the root quadratically. 

Because of practical concerns such as problem size or the lack of an analytical Jacobian, the Newton solver that is implemented is only an approximation of the ideal solver in \eqref{eq:theory:newton}.
For example, if a routine to evaluate the Jacobian is not provided, a common approach is to approximate it using finite difference methods.
This leads to the iteration scheme:
\begin{equation}
\label{eq:theory:newton-fd}
x_{n+1} = x_n - J_{FD}^{-1}(x_n)F(x_n),
\end{equation}
which is identical in form to \eqref{eq:theory:newton}, but replaces the Jacobian with a finite difference approximation (i.e., $J = J_{FD} + \text{error}$).
In larger problems, explicitly forming $J$ may be prohibitive (e.g., $J$ is a large and dense matrix), in which case Newton-Krylov methods, which only {\em approximately} solve the linear system associated with the correction step are often used~\cite{kelley2003solving}.
The results in the iteration scheme 
\begin{equation}
\label{eq:theory:newton-krylov}
x_{n+1} = x_n - Q(x_n)H^{-1}(x_n)Q^T(x_n)F(x_n) = x_n - J_{Krylov}^+(x_n)F(x_n),
\end{equation}
where $Q$ is the orthogonal subspace chosen by the Krylov solver,  $H$ is the approximate of $J$ in that subspace, and $+$ denotes the Moore-Penrose pseudoinverse.
In many applications, \eqref{eq:theory:newton-fd} or \eqref{eq:theory:newton-krylov} perform well, but as we will demonstrate shortly, there can be mismatches between the ideal and implemented solvers and systems that negatively impact the performance.

Finally, another common modification of Newton's method is the inclusion of a line search, which serves to limit the step size taken by Newton's method. 
There are a number of different line search strategies available~\cite{nocedal2006numerical}, we represent them via the inclusion of the state-dependent coefficient $\alpha$, resulting in the evolution law 
\begin{equation}
\label{eq:theory:newton-linesearch}
x_{n+1} = x_n - \alpha(x_n)J^{-1}(x_n)F(x_n),
\end{equation}
where $\alpha(x_n) = 1$ is a state where the full Newton step is taken, and where smaller values result in more conservative steps.
There are many other Newton-like methods; for example, one may combine a backtracking line-search with a Krylov solver~\cite{kelley2003solving}, or use other approaches such as Levenberg-Marquardt~\cite{nocedal2006numerical}, where the $\alpha$-like term appears elsewhere.

Due to the vast number of possibilities, we cannot enumerate the evolution laws associated with all of the nonlinear solvers that we may encounter even for a specific application such as the dynamic simulation of electronic circuit (Section~\ref{sec:circuits}).
In general however, our approach requires that the nonlinear solver of interest can be written as an {\em autonomous, deterministic, discrete-time dynamical system}, which all of the solvers listed above are examples of.
Given this requirement, solving a nonlinear system of equations is equivalent to finding a fixed point of this dynamical system.
If the approaches that will be outlined shortly are used, the evolution law associated with the solver (e.g.,~\eqref{eq:theory:newton}, \eqref{eq:theory:newton-fd}, \eqref{eq:theory:newton-krylov}, and \eqref{eq:theory:newton-linesearch}) does not need to be known explicitly, but such an equation must exist if the results are to be meaningful.

\subsection{Characterizing Solver Performance Using Linear Stability Analysis}
\label{sec:dynamical-system:stability}

In this section, we extend the analogy of ``solvers as dynamical systems'' to show that the convergence rate of the solver can be characterized by the eigenvalues obtained when one performs {\em linear stability analysis} on the resulting fixed point.
As mentioned previously, performing a full nonlinear analysis of such a system is likely to be extremely challenging due to the potentially fractal nature of the basin of attraction, and beyond the scope of what would be feasible on a black box system.
On the other hand, linear stability analysis is computationally tractable to perform even in large systems, and yields an upper bound on the convergence rate of the solver.

In what follows, we assume the nonlinear solver has the form,
\begin{equation}
\label{eq:theory:solver-form}
x_{n+1} = x_n  - \alpha(x_n) \tilde{J}^{-1}(x_n)F(x_n),
\end{equation}
where $F$ is the residual, $\tilde J$ is the approximation of the Jacobian used by the solver, and $\alpha$ is the step size of any associated line search.
The main quantity of interest is the spectrum of the linearization (or an empirical linear model) of \eqref{eq:theory:solver-form}.
This set of eigenvalues will enable us to determine the stability of the fixed point and quantify the asymptotic rate of convergence.

For the ideal implementation of Newton's method with $\tilde J = J$, the eigenvalues of the linearization of \eqref{eq:theory:solver-form} are zero because Newton's method converges quadratically.
The crux of our approach is that the ``mismatches'' between the theoretical and implemented numerical method can be traced to differences between the true and implemented Jacobians.
For the purpose of argument, we define the true Jacobian as $J=\frac{\partial F}{\partial x}$ even though it is certainly possible the routine that implements $F$ contains the errors.
We also assume the implemented Jacobian has the form: 
\begin{equation}
\label{eq:theory:jacobian-decomposition}
\tilde J  = J + \delta J + UCV,
\end{equation}
where $\delta J$ is a matrix whose entries are small, and the matrix-product $UCV$ is a low-rank matrix whose entries may be order 1.
This decomposition is somewhat arbitrary as terms can either be associated with $\delta J$ or $UCV$; all we intend to convey in \eqref{eq:theory:jacobian-decomposition} is that the Jacobian used in practice can deviate from the true Jacobian  in two ways -- there can be ubiquitous but small errors ($\delta J$) along with larger more serious deviations ($UCV$).  

To show the impact of each of these terms, we assume the fixed point is $x=0$, and compute the Taylor expansion of \eqref{eq:theory:solver-form}:
\begin{align}
\delta x_{n+1} &= I\delta x_n - \left(\alpha_0 + \nabla\alpha\cdot \delta x_n + \ldots\right)(J + \delta J + \hat J)^{-1}(J\delta x_n + \ldots)  \notag\\
&=   I\delta x_n - \left(\alpha_0 + \nabla\alpha\cdot \delta x_n + \ldots\right)\left((J + \hat J)^{-1}J\delta x_n - \delta JJ\delta x_n\right) + \mathcal{O}(\delta x^2) \notag\\
&= I \delta x_n - \alpha_0\left(J^{-1} - J^{-1}U\left(C^{-1} + VJ^{-1}U\right)^{-1}VJ^{-1}\right) J\delta x_n + \alpha_0 (\delta J)J\delta x_n + \mathcal{O}(\delta x^2)\notag\\ 
&=\left[(1-\alpha_0) I + \alpha_0 J^{-1}U\left(C^{-1} + VJ^{-1}U\right)^{-1}V +\alpha_0 (\delta J)J \right]\delta x_n + \mathcal{O}(\delta x^2).
\label{eq:theory:linearized-newton}
\end{align} 
For optimal performance, one would set $\tilde J = J$ and $\alpha_0 = 1$, which implies that $\delta x_{n+1} = \mathcal{O}(\delta^2)$, and that all of the eigenvalues of the linearized system are zero.
However, the presence of either $\delta J$ or $UCV$ will produce linear terms in the evolution law.
Although the eigenvalues themselves are both solver and system dependent, we can say something about their distribution based on the form of the system linearization in \eqref{eq:theory:linearized-newton}.
In particular, we can consider  the qualitative contribution of the three terms in the sum:
\begin{enumerate}
	\item The $(1-\alpha_0)I$ term, which is associated with the design of the line search, and shifts the system eigenvalues;
	\item The $\alpha_0 J^{-1}U\left(C^{-1}+ VJ^{-1}U\right)^{-1}V$ term, which is low-rank but can produce large shifts in the eigenvalues because the entries of $UCV$ are large;
	\item The $ \alpha_0 (\delta J)J$ term, which is full rank but should only perturb the eigenvalues  as the entries of $\delta J$ are assumed to be small.
\end{enumerate}

As a result, the spectrum of the linearized operator associated with a Newton-like method should have two parts: a ``cloud'' of small eigenvalues due to the 3rd term along with a few possible outliers due to the 2nd term.
For the purposes of numerical diagnostics, the cloud of eigenvalues represent the typical deviation of the implemented solver from the baseline.
For Newton's method, the baseline is that all the eigenvalues are at the origin, but damping or other modifications can change it.
The outlier eigenvalues are associated with what we refer to as {\em convergence anomalies}, which are the relatively small number eigenvalues that limit the convergence rate of the solver.

\subsubsection{Example: Rosenbrock Function}

To demonstrate a convergence  anomaly, we consider the problem of finding a minimum of the Rosenbrock function in $\mathbb{R}^{100}$.
The optimization problem is
\begin{equation}
\min_{x} \sum_{m=1}^{99} \left[100(x^{(m+1)}-(x^{(m)})^2)^2 + (x^{(m)} - 1)^2\right],
\end{equation}
which we will solve by finding an $x$ that satisfies the KKT conditions using Newton's method. 
The gradient of the objective function is 
{\small
\begin{equation}
\label{eq:theory:objective}
F^{(m)}(x) =
\begin{cases}
-400 (x^{(2)} - (x^{(1)})^2) x_1 + 2(x^{(1)} - 1) & \qquad m=1, \\
-400 (x^{(m+1)} - (x^{(m)})^2) x_m + 200 (x^{(m)} - (x^{(m-1)})^2)x^{(m)} + 2(x^{(m)} - 1) & \qquad 2\leq m \leq 99, \\
200 (x^{(100)} - (x^{(99)})^2)x^{(99)} & \qquad m=100.
\end{cases}
\end{equation}
}
Note that the global minimum of $F$ is given by the vector of 1s, and the Jacobian of $F$ is a tridiagonal matrix whose entries will be computed analytically and using a finite-difference method.

\begin{figure}[tp!]
	\centering 
	\includegraphics[width=0.7\textwidth]{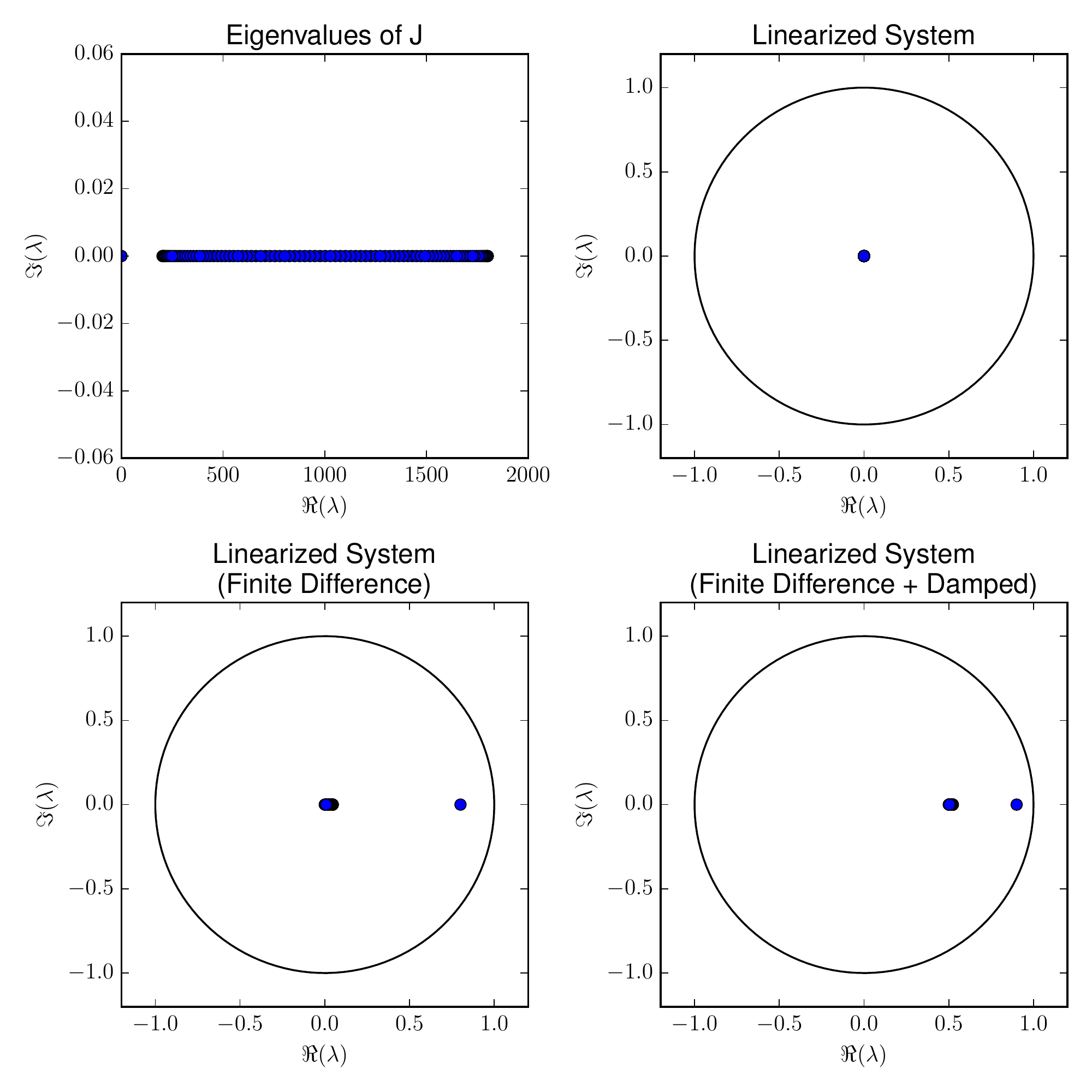}
	\caption{The top left plot shows the eigenvalues of the Jacobian for the 100-dimensional Rosenbrock equation about $x=1$.
		In addition to the clearly visible band of eigenvalues, there is a single isolated eigenvalue near $\lambda = 0.5$.
		The other plots show the eigenvalues associated with the linearization of Newton's method about $x=1$.
		The top right plot uses an analytical Jacobian with $\alpha=1$, the bottom left plot uses a finite difference approximation with $\Delta x = 0.01$ and $\alpha=1$, and the bottom right plot uses the same finite difference approximation with $\alpha = 0.5$.
		In all cases, the spectra consists of a cluster of eigenvalues, and, as demonstrated by the bottom row, a convergence  anomaly indicated by the exterior eigenvalue.
	}
	\label{fig:theory:rosenbrock}
\end{figure}

To demonstrate the impact that variants of Newton's method have on the spectrum of the system linearization for this problem, we implement  \eqref{eq:theory:objective} with a variety of Jacobian approximations.
The results of this procedure are shown in Fig.~\ref{fig:theory:rosenbrock}.
The top left plot in the figure shows the eigenvalues of $J$ (the Hessian of the objective function).
The smallest eigenvalue is a single outlier at $\lambda \approx 0.5$, but the majority of the eigenvalues lie in a band of $\lambda\in [200, 1850]$.
This system can be considered an example of a multi-scale system; in a neighborhood about $x=1$, the objective function increases rapidly in all directions except for the one corresponding to the eigenvector with the smallest eigenvalue. 
However, the scale separation is not large enough to cause numerical issues, and full Newton with an analytical Jacobian will converge rapidly given a good enough initial guess.
This is highlighted by the cluster of eigenvalues near zero ($|\lambda| < 10^{-5}$), and shown in the top right plot of Fig.~\ref{fig:theory:rosenbrock}.

The bottom row of Fig.~\ref{fig:theory:rosenbrock} shows the spectrum of the linearized Newton scheme when $J$ is approximated using a forward difference scheme with $\Delta x= 0.01$.
In both cases, there is a clear convergence  anomaly that is associated with the right-most eigenvalue in both plots.
The choice of $\Delta x = 0.01$ is larger than is reasonable, and was chosen so that the anomaly is clearly visible.
More realistic values such as $\Delta x = 10^{-4}$ still have a convergence  anomaly, but it is located  closer to the origin.
Less easy to see but equally important is the cluster of eigenvalues found near the origin in the bottom left plot or $\lambda\approx 0.5$ in the bottom right.
The radius of this cluster determines the typical deviation from the baseline, which for ``standard Newton'' is $\lambda=0$ but is $\lambda=0.5$ for damped Newton with $\alpha = 0.5$.

\subsection{Localizing Convergence Rate Anomalies}
\label{sec:theory:localize}

Once an anomaly has been detected, the corresponding eigenvectors can help us identify the source of the anomaly.
Without any further post-processing, these eigenvectors denote the states that converge the slowest to the fixed point, and the directions where errors will persist for the most iterations.
As a result, they help to pinpoint the system states that are {\em affected by the root cause of the anomaly}. 
However, if we have the ability to perform additional computational experiments, then it is possible that the cause of the anomaly can be further localized.

Ultimately, a convergence  anomaly is caused by a mismatch between the Jacobian of $F$ and the implemented Jacobian $\tilde J$ as shown in \eqref{eq:theory:jacobian-decomposition}.
For a small system, we can directly compare these two by approximating $J$ via a high-order finite difference method.
However, in the applications of interest to us, the state dimension can easily be large enough that explicitly computing a full approximation of $J$ is computationally intractable.
What we will show is that the eigenvectors associated with the isolated eigenvalues define an effective set of search directions that allow us to identify the image of the $U$ matrix, whose non-zero rows indicate the residual equations that are causing the convergence  anomaly (or equivalently, the residual equations that limit the convergence rate lie in the image of $U$).
While additional computational experiments/comparisons must be made, the number of required experiments is equal to the number of anomalies rather than the state dimension.

For the sake of simplicity, let $\alpha_0 = 1$ and $\delta J = 0$.
The inclusion of $\alpha_0 \neq 1$ shifts the eigenvalues but does not affect the eigenvectors, and hence will not impact the analysis.
The inclusion of $\delta J$ will perturb  the eigenvector of interest, but provided $\|\delta J\|$ is small and our assumption that this is truly an isolated eigenvalue holds, the effects of the perturbation on the resulting eigenvector will be small. 
As a result of the $\delta J$ term, there will be some ``noise'' in the prediction in practice that is being neglected here.

First, we want to establish the connection between the eigenvectors associated with convergence anomalies and the $UCV$ matrix.
With the additional simplifications, the leading order evolution law in some neighborhood of the root is:
\begin{equation}
x_{n+1} =  J^{-1}U\left(C^{-1} + VJ^{-1}U\right)^{-1}V x_n.
\end{equation}
As a result, the eigenvector, $v$, associated with an isolated eigenvalue, $\lambda$, satisfies
\begin{equation}
J^{-1}U\left(C^{-1} + VJ^{-1}U\right)^{-1}(Vv) = \lambda v.
\label{eq:theory:incorrect-eigenvector}
\end{equation}
The vector $v$ does not directly identify the equations where errors are present because it does not necessarily lie in the image of $U$ due to the $J^{-1}$ term.

To remove this term, we approximate the directional derivative of $F$ using a finite difference method in the direction of the eigenvector $v$.
Note that if a finite-difference method was used to construct $\tilde J$, the Jacobian used in the computation, then a higher order scheme should be used in this computation as our objective is to approximate the true Jacobian-vector product $J$.
Using \eqref{eq:theory:incorrect-eigenvector}, the directional derivative in the $v$ direction can be written as:
\begin{equation}
J v = \frac{1}{\lambda} J\left(J^{-1}U\left(C^{-1} + VJ^{-1}U\right)^{-1}V\right)v = \frac{1}{\lambda}U\left(C^{-1} + VJ^{-1}U\right)^{-1}Vv.
\label{eq:theory:identifying-equation}
\end{equation}
Once again, the exact value of the directional derivative is problem and error dependent, but it must lie in the image of $U$.
The only non-zero rows of $U$ are those associated with the equations causing the convergence anomalies.
As a result, if we can compute this matrix-vector product accurately, then the non-zero entries of $Jv$ indicate the equations causing the observed anomalies.
We should be clear that this approach will not identify or localize all errors, and it is certainly possible that the rate limiting error affects {\em all the residuals} or that the addition of ``noise'' due to the neglected $\delta J$ term masks some of the affected equations.
However, if we do detect a convergence limiting eigenvalue/eigenvector pair (or set of pairs) then  \eqref{eq:theory:identifying-equation} yields a single computational experiment that will identify the residuals affected by the dominant error(s).

\subsection{Tracking Convergence Anomalies in Homotopy Problems}

\begin{figure}[tp!]
	\centering
	\includegraphics[width=0.9\textwidth]{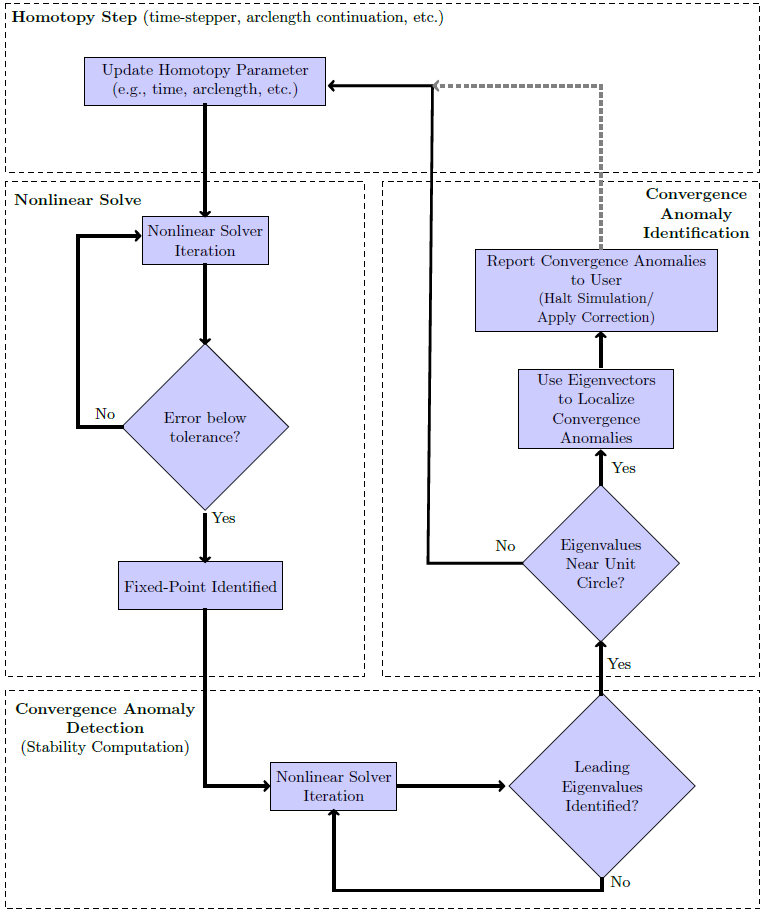}
	\caption{Flow chart of the ``additional experiments'' approach to convergence anomaly detection and identification.
		After every successful solve, additional solver iterations are performed to identify the leading eigenvalues and eigenvectors.
		If convergence anomalies, which appear as eigenvalues near the unit circle, are identified, then additional steps are taken to identify the equations/components responsible for the reduced convergence rate.}
	\label{fig:theory:flowchart-arnoldi}
\end{figure}

In the previous sections, we demonstrated that the performance of Newton's method applied to a single problem could be characterized by the eigenvalues of the system linearization.
However to obtain these eigenvalues, the fixed point must already have been identified.  
For problems with a stable fixed point, this analysis is applicable as a post-processing procedure that can characterize solver performance once a solution has been obtained.
In the context of homotopy problems, this procedure of fixed-point identification followed by analysis of the eigenvalues of the linearized system effectively results in a numerical bifurcation study.
For these problems, our objective is to track the convergence anomalies (i.e., the largest eigenvalues) as the simulation progresses (i.e., as one moves along the solution branch) to identify regions where poor solver performance is observed. 
Of particular interest are bifurcations in the solution branch, which imply the original branch of solutions is no longer stable (or no longer exists).

One possible result of a bifurcation is that the solver fails to converge beyond a certain time or parameter value, for example if a saddle-node bifurcation is encountered and there are no other solutions present in the system.
More subtle are cases where a bifurcation causes multiple stable solution branches to appear, and as a result, the uniqueness of the solution is lost as we do not know ahead of time which of the branches will be selected.
Furthermore because the convergence criterion for nonlinear solvers is often based on the norm of the residual, it is also possible for the nonlinear solving process to terminate near an unstable fixed point if the transient passes through a region that satisfies this criteria.
This is problematic in practice as minor changes to the convergence criteria or even the initial guess, which intuitively should have little impact on solver performance, can result in solver failure.
As a result, not all bifurcations signal the onset of solver failure, but in order to maintain the confidence and robustness of the solution, they should be avoided.

In order to detect a bifurcation, we require the eigenvalues of the linearization of Newton's method (or another solver) about the fixed point  at every point along the branch.
Because we do not have the explicit evolution equation for the combination of the nonlinear solver and system, we will use the approaches employed by black-box numerical continuation methods in order to obtain approximations of the needed eigenvalues.
First, we will outline the ``additional experiments'' approach, which uses additional computational experiments to compute the leading eigenvalues and eigenvectors.
Next, we will outline a data-driven method that approximates those quantities, but can be applied to a larger class of problems.

\subsubsection{The Additional Experiments Approach to Numerical Diagnostics}
The additional experiments approach is shown in Fig.~\ref{fig:theory:flowchart-arnoldi}, and can be thought of as a parallel to the process that occurs in pseudoarclength continuation.
The procedure we propose consists of four steps:
\begin{enumerate}
	\item First is the homotopy step, which updates the parameter being swept.
	In time-integration problems, this is time, and in pseudoarclength continuation, this is the approximation of the arclength.
	\item
	At each homotopy step, a nonlinear system of equations must be solved, which for the solvers we consider, requires running a number of Newton-like iterations until the solution has been identified.
	This solution is the fixed point of the iterative map associated with our nonlinear solver.
	\item In order to identify the convergence anomalies, we require the leading eigenvalues (i.e., those largest in magnitude).
	To do this, we run an additional set of computational experiments by systematically perturbing the system away from the identified fixed point. 
	Because we are only interested in the leading eigenvalues, these perturbations can be chosen using the Implicitly Restarted Arnoldi Method (IRAM), which is implemented in ARPACK~\cite{lehoucq1998arpack}, though we currently compute an approximation of the entire Jacobian using a finite difference scheme.
	\item If any of the eigenvalues are near the unit circle, then additional experiments can be performed to localize these convergence anomalies to a subset of the residual equations.
	After this point, the user can terminate the procedure, modify the system definition, apply corrections by adjusting solving process, or simply passively track the anomalies. 
\end{enumerate}
In the examples that will follow, we choose to passively observe the effects of the errors and other convergence anomalies even if they will ultimately result in solver failure. 
In practice, however, one might want to apply some form of corrective action in order to prevent solver failure/mitigate the effects of the convergence  anomaly.

To use this approach, we must be able to run additional iterations of the nonlinear solver starting from pre-specified initial conditions.
The benefit of performing additional experiments is that highly accurate approximations of the eigenvalues and eigenvectors can be computed regardless of how many iterations the nonlinear solver required to find the solution.
However, not all implementations, particularly those in commercial software packages, enable this.
When additional experiments are not possible, the data-driven approach that will be described should be used instead.

\subsubsection{A Data-Driven Approach to Numerical Diagnostics}

\begin{figure}[tp!]
	\centering
	\includegraphics[width=0.9\textwidth]{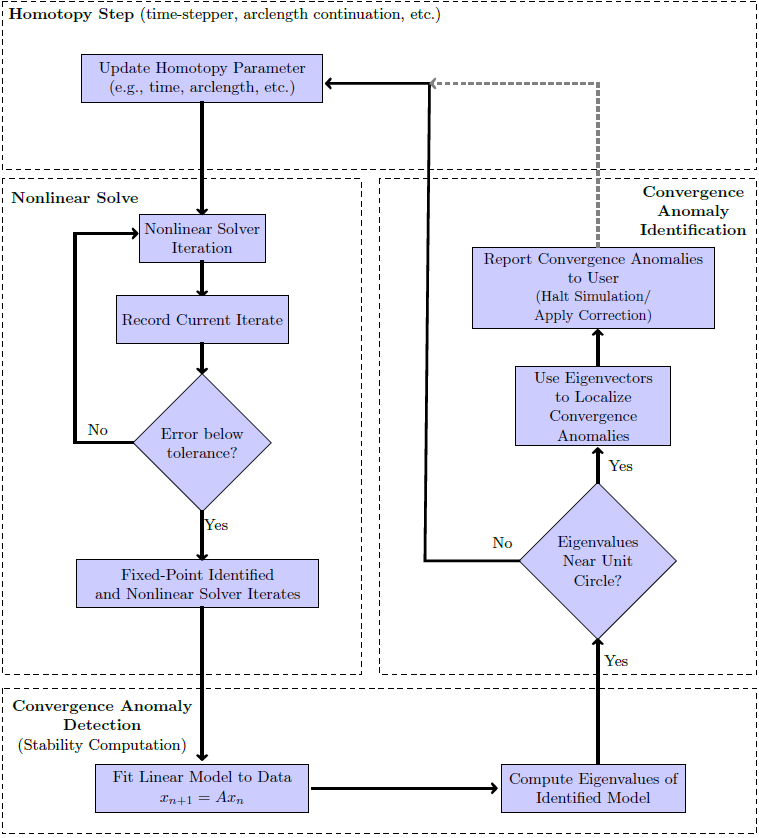}
	\caption{Flow chart of the data-driven approach to convergence anomaly detection and identification.
		After each internal iteration of the nonlinear solver, the current iterate is recorded.
		Once a fixed-point has been identified, these data are used to fit a linear model whose eigenvalues approximate that of the system linearization about the fixed point.
		If convergence anomalies, which appear as eigenvalues near the unit circle, are identified, then additional steps are taken to identify the equations/components responsible for the reduced convergence rate}
	\label{fig:theory:flowchart-data}
\end{figure}

In order to apply this technique to problems where additional computational experiments cannot be run, we also present a data-driven approach to obtain the same quantities.
This approach borrows techniques from data-driven system identification, and fits linear models to snapshots of the system state that are saved at every step of the nonlinear solve.
As shown in the flowchart in Fig.~\ref{fig:theory:flowchart-data},  this procedure consists of the same four steps as the ``additional experiments'' approach.
The differences between the two approaches are in steps 2 and 3.
In particular,
\begin{enumerate}
	\item[2.] At each homotopy step, a nonlinear system of equations must be solved. 
	After every internal iteration of the nonlinear solver, we save the current state (i.e., the current solution ``guess'' at each iterate).
	The data must be ordered by iteration number, as we will apply time-series methods to it in the subsequent step.
	Once the solver has converged, this data is passed along with the numerically identified solution to the third step of the procedure.
	\item[3.] In this step, our objective is to fit a linear model to the time series data.
	Because the data are not centered around zero, we first take the difference between each pair of snapshots:
	\begin{equation}
	\begin{bmatrix}
	x_1 & x_2 & x_3 & \cdots & x_M
	\end{bmatrix}
	\Rightarrow
	\begin{bmatrix}
	\Delta x_1 & \Delta x_2 & \Delta x_3 & \cdots & x_{M-1}
	\end{bmatrix}
	\text{ where } \Delta x_m = x_{m+1} - x_m,
	\end{equation}
	which reduces the size of the data set by one, but allows us to fit a linear rather than an affine model to the data, which we do using Dynamic Mode Decomposition (DMD)~\cite{schmid2010dynamic,tu2013dynamic}.
	In this context, we will assume our data are collected near the fixed point, but the connections between DMD and Koopman operator theory for nonlinear systems~\cite{tu2013dynamic,williams2014data} and the global stability analysis enabled by that framework~\cite{mauroy2013spectral} suggest this approach could still have mathematical meaning without this restriction on the data. 
	More concretely, we identify the matrix $A$ that solves the least squares problem:
	\begin{equation}
	\min_{A} \sum_{m=1}^{M-1} \|\Delta x_{m+1} - A \Delta x_m\|^2,
	\end{equation}
	while simultaneously minimizing the Frobenius norm of $A$.
	We then compute the eigenvalues and eigenvectors of $A$ using standard approaches, which are provided to the fourth step of the extra experiments approach for convergence anomaly identification.
\end{enumerate}
The data-driven approach can be applied to problems where the nonlinear solve process can be passively observed but not directly altered, which is often the case when the solution process is part of a larger software package. 
However, it requires data, so we must have the ability to obtain information at every iterate of the solving process.

Regardless of whether a data-driven or additional experiments approach is taken, the goal of the methods we proposed is to produce an approximation of the eigenvalues of the system linearization at each step of the homotopy procedure.
By tracking these eigenvalues we can determine when convergence anomalies appear and disappear, and characterize whether or not the current solution branch will remain stable for the region of interest or bifurcate.
It should be noted that other features that may be useful for numerical diagnostics, such as discontinuity induced bifurcations (DIBs)~\cite{di2008discontinuity, di2010discontinuity} due to loss of smoothness in the underlying dynamical system, are also indicated by the eigenvalues.
For example, border collisions are appear as sudden ``jumps'' in the eigenvalues of the system linearization.
As a result, the spectrum contains useful information that may help diagnose other types of numerical issues than the the ones we will focus on here.

\section{Diagnosing Numerical Failures in Electrical Circuit Simulation}
\label{sec:circuits}

To demonstrate the use of these methods, we apply them to circuit simulation applications.
This is an interesting application of such methods for two reasons: (1) techniques such modified nodal analysis naturally produce large systems of DAEs, which require nonlinear equations to be solved at each time step, and (2) when the system is assembled using an object-oriented component-composition approach, errors in the system implementation can be isolated to individual components not just the residual equation.
In the remainder of this section, we will give a brief outline of how the governing equations for circuit simulation can be assembled from individual components.
Then we will present three examples highlighting our approach's ability to identify and isolate implementation and system-level issues in electrical circuits.

\subsection{Circuit Simulation via Component Composition}

One method of simulating electrical circuits is through the use of {\em modified nodal analysis}~\cite{ho1975modified}, which provides an algorithmic way of constructing the governing equations of large circuits given component-level residuals and Jacobians.
Compared to the overall system, whose dynamics may be quite complex, the component-level residuals and Jacobians tend to be simpler to derive and compute.

Figure~\ref{fig:motivation:components} shows the component level residuals and Jacobians for a resistor, an ideal voltage source, a capacitor, and an ideal diode respectively.
Physically, the residuals associated with the resistor and diode are mappings from the voltage drop across the component to the current leaving the nodes indicated in the figure.
Both of these components can be parameterized: $R$, the resistance of the resistor, will take on a wide range of values, but the diode parameters will remain fixed at $I_S = 10^{-12}$ A, $n=1$, and $V_T = 26$ mV. 
Components such as the idealized voltage source can also have internal states and constraints, which in that example is $I_V$, the current produced by the voltage source.  
To compensate for these additional state variables, additional residual equations associated with these variables, such as the third equation in Fig.~\ref{fig:motivation:components}, are introduced; we refer to these new equations as internal nodes.
These residuals do not necessarily have an interpretation as a mapping from voltages to currents; in this case, it enforces that the voltage drop across the source is $V(t)$.
Finally, components such as capacitors depend upon the time-derivative of the nodal voltages, which are approximated in terms of state variables and prior solutions when methods like BDFs are used to solve the resulting ODE/DAE.

\begin{figure}[tp!]
	\centering

	\begin{tikzpicture}[every node/.style={outer sep=0.0mm, inner sep=0.0mm}]
	\node[label= left:{Node 0 ($V_0$)\hspace*{1em}}] (v2) at (0, 3) {};
	\node[label= right:{\hspace{1em}Node 1 ($V_1$)}] (v3) at (5, 3) {};
	\draw (v2) to[resistor=$R$, *-*] (v3);
	\end{tikzpicture}
	
	\begin{equation*}
	r_R(V_0, V_1)
	= 
	\begin{bmatrix}
	\frac{V_0 - V_1}{R} \\
	-\frac{V_0 - V_1}{R} \\
	\end{bmatrix}, \quad
	J_R(V_0, V_1)
	= 
	\begin{bmatrix}
	\frac{1}{R} & -\frac{1}{R}  \\
	-\frac{1}{R} & \frac{1}{R} 
	\end{bmatrix}
	\end{equation*}

	\begin{tikzpicture}[every node/.style={outer sep=0.0mm, inner sep=0.0mm}]
	\node[label= left:{Node 0 ($V_0$)\hspace*{1em}}] (v2) at (0, 3) {};
	\node[label= right:{\hspace{1em}Node 1 ($V_1$) }] (v3) at (5, 3) {};
	\node[label={Internal Node  ($I_V$)}] at (2.5, 2.0) {};
	\draw (v2) to[sV=$V(t)$, *-*] (v3);
	\end{tikzpicture}
	
	\begin{equation*}
	r_V(V_0, V_1, I_V)
	= 
	\begin{bmatrix}
	I_V \\
	-I_V \\
	V_0 - V_1 - V(t)
	\end{bmatrix}, \quad
	J_V(V_0, V_1)
	= 
	\begin{bmatrix}
	0 & 0 & 1  \\
	0 & 0 & -1 \\
	1 & -1 & 0 
	\end{bmatrix}
	\end{equation*} 
	
	\begin{tikzpicture}[every node/.style={outer sep=0.0mm, inner sep=0.0mm}]
	\node[label= left:{Node 0 ($V_0$)\hspace*{1em}}] (v2) at (0, 3) {};
	\node[label= right:{\hspace{1em}Node 1 ($V_1$) }] (v3) at (5, 3) {};
	\draw (v2) to[capacitor=$C$, *-*] (v3);
	\end{tikzpicture}
	
	\begin{equation*}
	r_C(V_0, V_1)
	= 
	\begin{bmatrix}
	C(\dot V_0 - \dot V_1)\\
	-	C(\dot V_0 - \dot V_1) 
	\end{bmatrix}\approx 
	\begin{bmatrix}
	\alpha C(V_0 -  V_1) - C(\beta_0 - \beta_1) \\
	-	\alpha C( V_0 -  V_1)  + C(\beta_0 - \beta_1)
	\end{bmatrix}
	, \quad
	J_C(V_0, V_1)
	= 
	\begin{bmatrix}
	\alpha C & -\alpha C  \\
	-\alpha C & \alpha C  
	\end{bmatrix}
	\end{equation*} 
	
	\begin{tikzpicture}[every node/.style={outer sep=0.0mm, inner sep=0.0mm}]
	\node[label= left:{Node 0 ($V_0$)\hspace*{1em}}] (v2) at (0, 3) {};
	\node[label= right:{\hspace{1em}Node 1 ($V_1$)}] (v3) at (5, 3) {};
	\draw (v2) to[diode, *-*] (v3);
	\end{tikzpicture}
	
	{\small
	\begin{equation*}
	r_D(V_0, V_1)
	= 
	\begin{bmatrix}
	I_S\left(\exp\left(\frac{V_0 - V_1}{nV_T}\right) - 1 \right) \\
	-I_S\left(\exp\left(\frac{V_0 - V_1}{nV_T}\right) - 1 \right)  \\
	\end{bmatrix}, \quad
	J_D(V_0, V_1)
	= 
	\begin{bmatrix}
	\frac{I_S}{nV_D} \exp\left(\frac{V_0 - V_1}{nV_T}\right)& -\frac{I_S}{nV_T} \exp\left(\frac{V_0 - V_1}{nV_D}\right) \\
	-\frac{I_S}{nV_D} \exp\left(\frac{V_0 - V_1}{nV_T}\right) &\frac{I_S}{nV_T} \exp\left(\frac{V_0 - V_1}{nV_T}\right) 
	\end{bmatrix}
	\end{equation*}
	} 
	
	\caption{Component-level residuals and Jacobians (with respect to the states) for a resistor, an ideal voltage source, a capacitor, and an ideal diode.
	In the equations above, $R$ is the resistance, $V(t)$ is the voltage of the source, $C$ the capacitance of the capacitor, $I_S$ is the saturation current of the diode, $n$ is the emission coefficient, and $V_T$ is the thermal voltage.
	For the capacitor, the time derivatives of voltage are approximated as $\dot V_i = \alpha V_i - \beta_i$, which is consistent with the approach taken when BDFs are used to solve ODEs/DAEs.
	Physically, the residual equations either map the voltage drop across a component (or the time derivative of that drop) to the current through it, or as in the third equation of the voltage source, denote constraints that must be satisfied.
	Using modified nodal analysis, these individual components can be aggregated to determine the governing equations of a larger circuit.
	}
	\label{fig:motivation:components}
\end{figure}

This component-level information is then aggregated into a single system-level set of residuals and the associated Jacobian.
The system state is defined by the voltage at a set of pre-defined nodes augmented by the internal states of the components.
Using these voltages along with prior knowledge about which components connect the nodes, the contribution of each component level residual and Jacobian can be mapped to the appropriate system level residual and Jacobian entries.
As a result, we assume the system-level residual, $F$, has the form
\begin{equation}
F_m(x) = \sum_{n=1}^N a_{mn} r_n(x),
\label{eq:composition:residual}
\end{equation}
where $N$ is the number of component-level residuals, $F_m$ is the $m$th system-level residual, $a_{mn}$ is either 0 or 1 and governs the connection of components to nodes, and $r_n$ is the $n$th component-level residual equation.
Then, the Jacobian of $F$ can be written as 
\begin{equation}
J = A\frac{\partial r}{\partial x} = AR,
\label{eq:composition:jacobian}
\end{equation}
where $A$ contains the $a_{mn}$ coefficients, and $R$ contains the partial derivatives of the component-level residual equations. 
Note that $A$ is almost certainly a rectangular matrix as we typically have many more components than nodes.

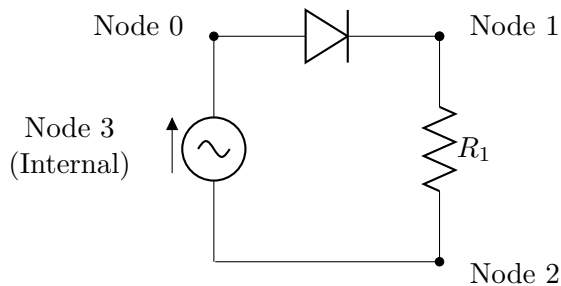
\begin{figure}
	\centering
	\usetikzlibrary{circuits.logic.US,circuits.logic.IEC}
	\begin{tikzpicture}[every node/.style={outer sep=0.0mm, inner sep=0.0mm}]
	
	\node (v1) at (0, 0) {};
	\node[label=above left:{Node 0\hspace*{1em}}] (v2) at (0, 3) {};
	\node[label=above right:{\hspace{1em}Node 1}] (v3) at (3, 3) {};
	\node[label=below right:{\hspace{1em}Node 2}] (v4) at (3, 0) {};
	
	\draw (v1) to[sV={\parbox{1in}{\centering Node 3 \\ (Internal)}}] (v2);
	\draw (v2) to[diode=$D_1$, *-*] (v3);
	\draw (v3) to[resistor=$R_1$, *-*] (v4);
	\draw (v4) -- (v1);
	\end{tikzpicture}
	\caption{Simple electrical circuit consisting of a diode and resistor in series with a time-varying voltage source whose time-evolution can be computed by solving the nonlinear system of equations in \eqref{eq:motivation:simple-circuit-residual}.
	}
	\label{fig:motivation:circuit}
\end{figure}

As an example of this approach, consider the circuit in Fig.~\ref{fig:motivation:circuit}, which consists of a time-varying voltage source connected to a diode and resistor in series. 
Using the component-level residuals, the system level residual can be written as:
{\small
\begin{equation}
F(V_0, V_1, V_2, I_V) = 
\begin{bmatrix}
1 & 0 & 0 & 0 & 1 & 0 & 0 \\
0 & 1 & 1 & 0 & 0 & 0 & 0 \\
0 & 0 & 0 & 1 & 0 & 1 & 0 \\
0 & 0 & 0 & 0 & 0 & 0 & 1 \\
\end{bmatrix}
\begin{bmatrix}
r_D(V_0, V_1) \\
r_R(V_1, V_2) \\
r_V(V_0, V_1, I_V)
\end{bmatrix}
= 
\begin{bmatrix}
I_S\left(\exp\left(\frac{V_0 - V_1}{nV_T}\right) - 1 \right)  + I_V \\
-I_S\left(\exp\left(\frac{V_0 - V_1}{nV_T}\right) - 1 \right) + \frac{V_1 - V_2}{R_1} \\
-\frac{V_1 - V_2}{R_1} - I_V \\
V_0 - V_2 - V(t)
\end{bmatrix},
\label{eq:motivation:simple-circuit-residual}
\end{equation}
}
and the voltages/currents obtained by solving $F(V_0, V_1, V_2, I_V) = 0$
\footnote{Because the circuit  in Fig.~\ref{fig:motivation:circuit} lacks a ground, $F=0$ does not have a unique solution.
In practice, one of the three voltages $V_0$, $V_1$, or $V_2$ are set to zero in order to set the reference level, and the corresponding equations and Jacobian entries removed.
}.
Similarly, the system level Jacobian can be written as 
\begin{align}
J(V_0, V_1, V_2, I_V) &= 
\begin{bmatrix}
1 & 0 & 0 & 0 & 1 & 0 & 0 \\
0 & 1 & 1 & 0 & 0 & 0 & 0 \\
0 & 0 & 0 & 1 & 0 & 1 & 0 \\
0 & 0 & 0 & 0 & 0 & 0 & 1 \\
\end{bmatrix}
\begin{bmatrix}
\hat J_D(V_0, V_1) \\
\hat J_R(V_1, V_2) \\
\hat J_V(V_0, V_1, I_V)
\end{bmatrix} \\
&= 
\begin{bmatrix}
\frac{I_S}{nV_D} \exp\left(\frac{V_0 - V_1}{nV_T}\right) & -\frac{I_S}{nV_D} \exp\left(\frac{V_0 - V_1}{nV_T}\right) & 0 & 1 \\
-\frac{I_S}{nV_D} \exp\left(\frac{V_0 - V_1}{nV_T}\right) & \frac{I_S}{nV_D} \exp\left(\frac{V_0 - V_1}{nV_T}\right) + \frac{1}{R_1} & -\frac{1}{R_1} & 0 \\
0 & -\frac{1}{R_1} & \frac{1}{R_1} & -1 \\
1 & 0 & -1 & 0
\end{bmatrix}
,\notag
\label{eq:motivation:simple-circuit-jacobian}
\end{align}
where $\hat J_R$, $\hat J_D$, and $\hat J_V$ are the component-level Jacobians from Fig.~\ref{fig:motivation:components} ``padded'' with zeros to account for nodes that are not associated with the component.
In practice, one evaluates $J_R$, $J_D$, and $J_V$ and then distributes the elements rather than padding the matrix, but the end result is the same as what is shown above.
These governing equations could be assembled by hand, and the main benefit of this approach is that they can be algorithmically assembled and scales well to larger number of nodes and components.
This also allows us direct access to the vector of $r_n$ and the component-level Jacobians aggregated in $R$, which we will use to further localize errors down to individual components.

Identifying component-level issues is done using the same approach outlined in Sec.~\ref{sec:theory:localize}, but applied to component rather than system-level Jacobians and residuals.
The system level residual can be written in the form shown in \eqref{eq:composition:residual}, and the ``ideal'' Jacobian is of the form shown in \eqref{eq:composition:jacobian}.
The various approximations used to construct $J$ in Sec.~\ref{sec:theory:localize} are still used in this context, but are now applied at a component level rather than at the system level.
As such, we write 
\begin{equation}
\tilde J = J + \delta J + \hat{J} = A(R + \delta R +  U_rC_rV_r) = A\tilde{R},
\end{equation}
where $U_rC_rV_r$ contain the components residuals with error, the magnitudes of those errors, and the nodes that contribute respectively.
Assuming once again that the $\delta R$ term is negligible, our objective is to identify the image of $U_r$, which contains the components residual equations where these rate-limiting errors occur.
To do this, we must have access to the code that computes $\tilde{R}$, the Jacobian of the component-level residuals, and the code that evaluates the component-level residuals themselves.
Then we compare the directional derivative obtained from $\tilde{R}$ with the directional derivative computed from $r$ using a high order scheme:
\begin{equation}
\frac{r(x + \epsilon v) - r(x-\epsilon v)}{2\epsilon} - \tilde{R}v =- U_r C_r V_rv    - \delta R v + \mathcal{O}(\epsilon^2),
\end{equation}
where the $\epsilon^2$ is due to our choice of finite difference method and could be adjusted by choosing a difference finite difference scheme.
The resulting vector lies in the span of $U_r$, which signifies the component-level equations where issues are present, with the addition of ``noise'' due to the $\delta R$ and $\mathcal{O}(\epsilon^2)$ terms.
There are two additional perturbations, the inclusion of the $\delta R v$ term and the error associated with the finite difference scheme that will introduce noise into the resulting computation that may mask equations with small errors.
Furthermore, we often normalize by the value of $r$ to account for scale differences in the residual equations.
With these preprocessing steps in place, any dominant peaks in the resulting vector will indicate the components in the range of $U_r$, which is the equations where our rate-limiting errors occur.

In the end, composed systems like the circuits we will consider can be treated using the same numerical techniques discussed in Sec.~\ref{sec:dynamical-systems}.
However, because we know something about the form of the system, in this case that is the linear combination of component-level residuals and Jacobians, we can exploit this knowledge to further localize the cause of a convergence  anomaly.
The overall procedure remains the same: eigenvalues are used to detect convergence issues and the corresponding eigenvector narrows the issue down to the affected nodes, but now we can perform an additional step to isolate issues to the level of individual components.

\subsection{Example: Implementation Error Localization in a Diode Bridge}

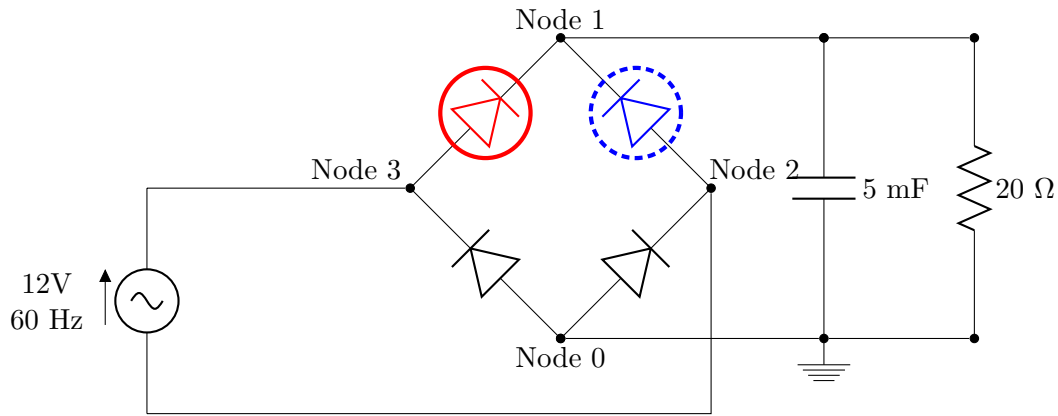
\begin{figure}[tp]
	\centering
	\usetikzlibrary{circuits.logic.US,circuits.logic.IEC}
	\begin{tikzpicture}[every node/.style={outer sep=0.5mm, inner sep=0.0mm}]
	
	\node[label=above left:{Node 3}] (v1) at (3.5,3) {};
	\node[label={Node 1}] (n2) at (5.5,5) {};
	\node[label=below:{Node 0}](n3) at (5.5,1) {};
	\node[label=above:{$\quad\qquad$Node 2}] (v2) at (7.5,3) {};
	\draw (v2) -- (7.5,0)  -- (0,0)  to[sV={\parbox{0.5in}{\centering 12V \\ 60 Hz}}] (0, 3) -- (v1);
	
	\draw (n3) to[diode, *-*] (v2);
	\draw (v1) to[diode, color=red] (n2);
	\draw[ultra thick, color=red] (4.5, 4) circle (0.6cm);
	\draw (v2) to[diode, color=blue] (n2);
	\draw[ultra thick, densely dashed, color=blue] (6.5, 4) circle (0.6cm);
	\draw (n3) to[diode, *-*] (v1);
	\node (v3) at (9,5) {};
	
	\node (v4) at (11,5) {};
	
	\draw (n2) to[short, *-*] (v3);
	\draw (v3) -- (v4);
	\draw[fill=none] (11,5) to[resistor={20 $\Omega$}, *-*] (11,1) node (v6) {};
	\draw (v3) to[capacitor={5 mF}, *-*] (9,1) node (v5) {};
	\draw (n3) -- (v5);
	\draw (v5) -- (v6);
	
	\draw (v5) node[ground]{};
	\end{tikzpicture}
	
	\caption{Circuit diagram of the diode bridge. 
		Two of the diodes are colored in red and blue (solid and dashed circles respectively), and will have implementation errors introduced into their Jacobians in the examples that follow.
		The each of the diodes are implemented with a parallel resistor, which are not shown, in order to include a minimum conductance. 
		These extra resistors will be treated like independent components, and therefore have separate component-level residuals.
	}
	\label{fig:circuit:bridge-error}
\end{figure}

In the next two examples, we apply our approach to the diode bridge shown in Fig.~\ref{fig:circuit:bridge-error}.
The diode bridge system consists of {\em five states}, the voltages of nodes 0-3 as labeled in Fig.~\ref{fig:circuit:bridge-error}, and a node that is internal to the voltage source
The components in the system consists of a 12 V, 60 Hz ideal sinusoidal voltage source, a diode bridge consisting of the four diodes and the 0.005 F filter capacitor, and a 20 $\Omega$ resistive load. 
To aid in convergence, simulators like SPICE~\cite{Nagel:M382} sometimes implement minimum conductances on nonlinear devices, which appear as large parallel resistors, and we do the same here by placing a $10^{12}$ $\Omega$ resistor in parallel with all diodes unless otherwise specified.
Because of the capacitor, the system is a differential algebraic equation (DAE) rather than an algebraic system of equations (see Fig.~\ref{fig:motivation:components}), and will be integrated over the 20 ms time window of interest using the 1st order Backwards Differentiation Formula (BDF)~\cite{shampine1979user}.
The nonlinear solve required at each time step is performed using the Newton solver implemented by NOX~\cite{heroux2005overview} with an analytical Jacobian.

As a benchmark, we will first present the diode bridge system without deliberately introduced errors.
Next, we will include sign flips in the Jacobians of the resistors introduced to produce the minimum conductance in the blue and red diodes.  
This combination of errors will result in solver divergence in the 20 ms window of interest.
Next, we will correct the implementation error associated with the blue diode, but retain the error in the red diode.
In this case, the simulation successfully converges over the 20 ms window, but there are numerically significant deviations from the true solution.
In our last example, we focus on the case where no implementation errors are present, but the resistors that implement the minimum conductance of the diodes are removed.
Our objective in all examples is to demonstrate how the approach outlined in Section~\ref{sec:dynamical-systems} signals these issues may be occurring, and show it can be used to localize the sources of error down to individual components when the convergence  anomaly is due to a subset of the components.

\begin{figure}[t!p]
	\centering
	\includegraphics[width=0.9\textwidth]{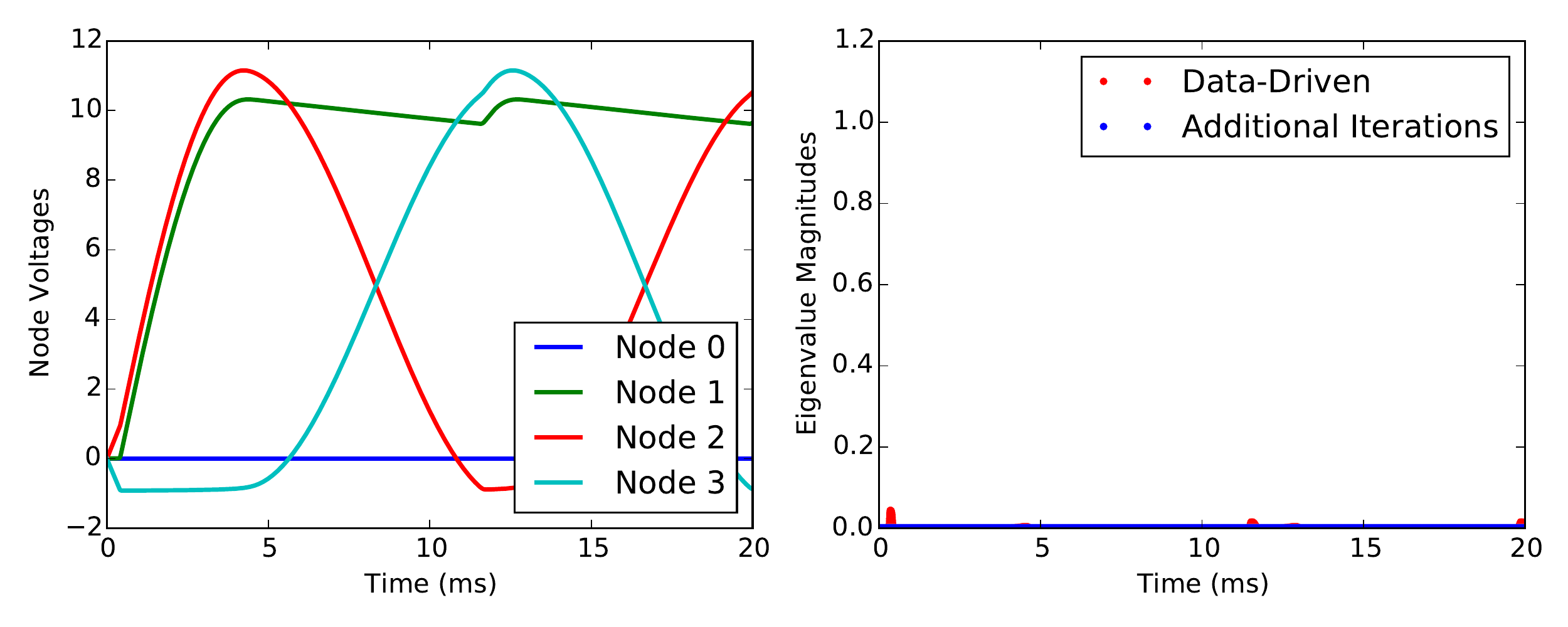}
	\caption{(left) The voltages at the labeled nodes over the 20 ms time window.
		(right) The eigenvalues obtained using the ``additional experiments'' approach and the ``data-driven'' approach for the nonlinear solve associated at each time step.
		Due to many Newton-solves terminating after only a single iteration, the data-driven approach, which is denoted by the red markers, can only be applied to a minority of the time steps.
		In both examples, the eigenvalues remain clustered near the origin.
	}
	\label{fig:circuit:bridge-fine}
\end{figure}

\subsubsection{Reference Case}
Figure~\ref{fig:circuit:bridge-fine} shows the evolution of the node voltages over the 20 ms window of interest in the absence of implementation errors.
Nodes 0, 2, and 3 have predictable behaviors based on the system configuration: node 0 is tied to ground at all times, and the potential difference between nodes 2 and 3 must be a 60 Hz sinusoid due to the ideal voltage source.
The voltage across the load, which is determined by node 1, is ``DC-like'' with ripples whose magnitude is determined by the load and the capacitance of the filter capacitor~\cite{horowitz1989art}. 

The right plot of Fig.~\ref{fig:circuit:bridge-fine} shows the magnitude eigenvalues of the linearized Newton method for each time step of the procedure.
The blue markers denote the eigenvalues obtained by the additional experiments approach, and the red markers show the eigenvalues obtained by the data-driven approach.
Because we need at least three Newton iterations to approximate the eigenvalues, the gaps in the figure indicate the times where fewer than three Newton iterations were performed.
However, when sufficient data exists to compute an approximation, both methods consistently generate a set of eigenvalues near the origin, which implies there are no convergence anomalies.

\subsubsection{Two Implementation Errors}
\begin{figure}[tp!]
	\centering
	\includegraphics[width=0.6\textwidth]{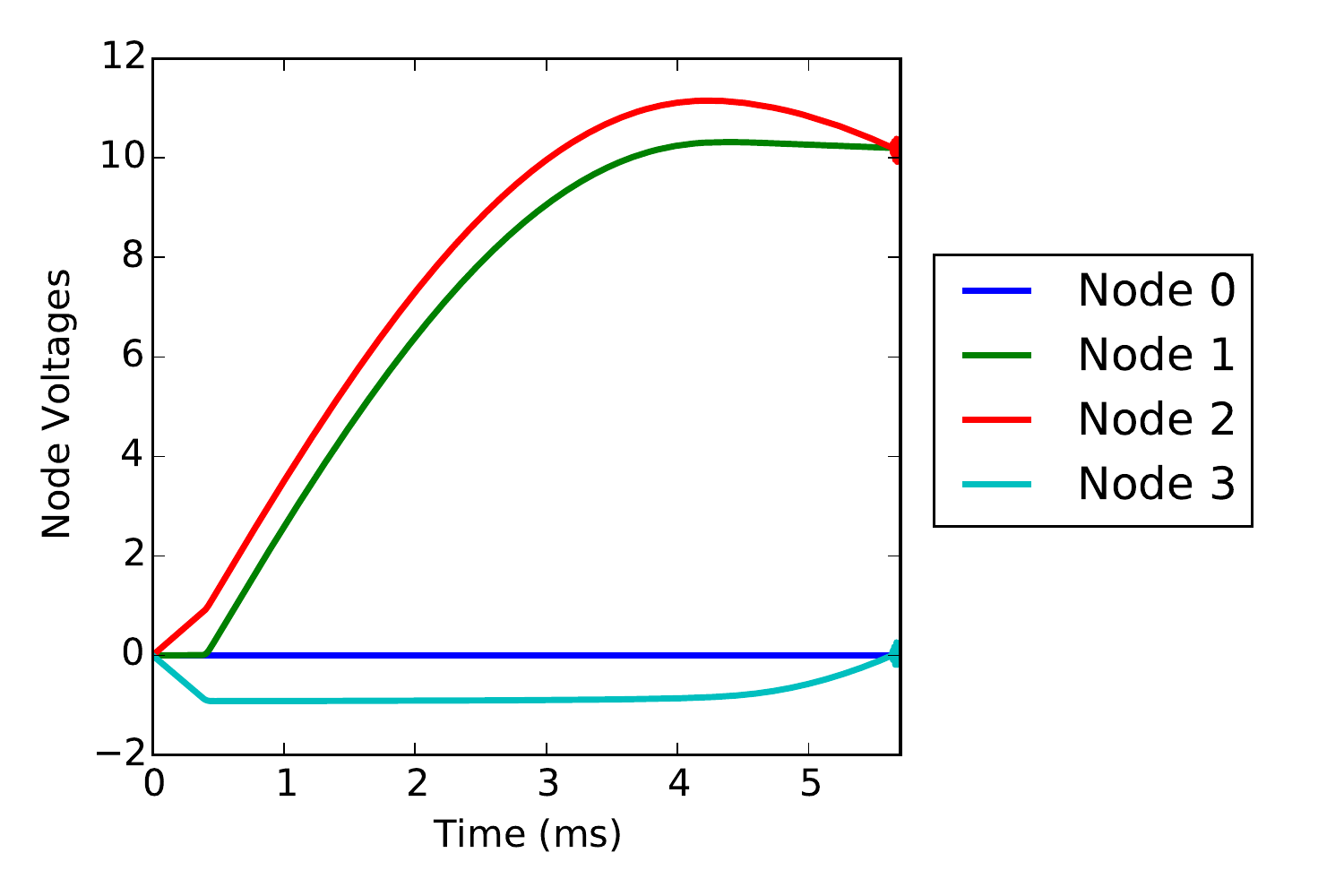}
	\caption{Plot of the node voltages as a function of time.
		After $t\approx 5.7$ ms, the Newton solver fails to converge and the time-integration procedure terminates.
	}
	\label{fig:circuit:bridge-failure}
\end{figure}

Now we consider the case where the diode bridge shown in Fig.~\ref{fig:circuit:bridge-error} has two components -- the red and blue diodes (solid and dashed circles respectively) -- with sign flips in their Jacobians, which are meant to be prototypical programming errors.
This combination of errors will result in complete solver failure (i.e., return ``NaNs'') in the 20 ms time-window desired.
The segment of the trajectory that was computed in shown in Fig.~\ref{fig:circuit:bridge-failure}.
Note the oscillations that appear just before failure at $t\approx 5.7$ ms.

\begin{figure}[tp!]
	\centering
	\includegraphics[width=0.9\textwidth]{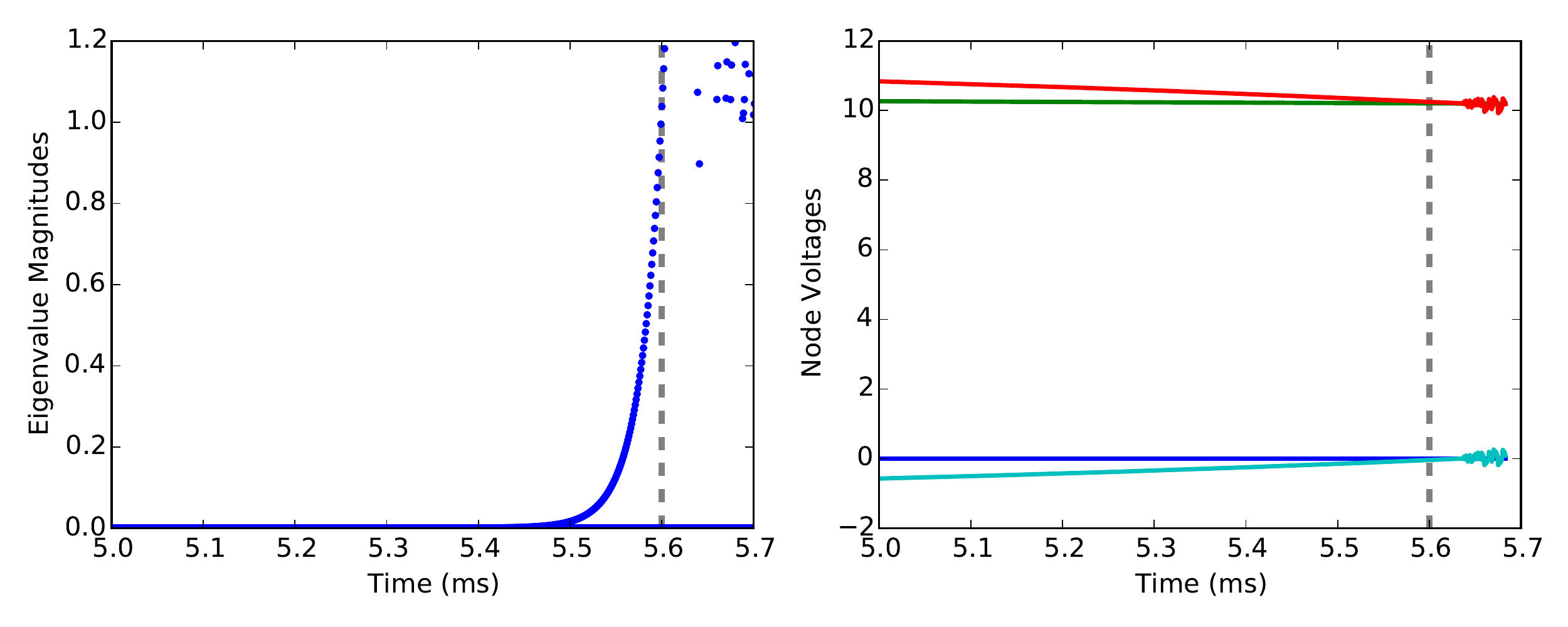}
	\caption{
		(left) The magnitude of the eigenvalue obtained using the additional experiments approach near the solver failure at 5.7 ms. 
		There is a period-doubling bifurcation near 5.6 ms that indicates the convergence rate of the solver has slowed, and that failure may be imminent.
		(right) The voltage evolution in the same time window.
		The gray line in both plots indicates the approximate location of the period-doubling bifurcation, which precedes the visible loss of convergence closer to 5.7 ms.
	}
	\label{fig:circuit:bridge-failure-comparison}
\end{figure}

At the current time, the approach  does not feedback into the underlying solver/system.
Our goal  is to  diagnose the convergence anomalies observed rather than to prevent failure.
As such, Fig.~\ref{fig:circuit:bridge-failure-comparison} shows the magnitude of the eigenvalues obtained using the additional experiments approach as a function of time (e.g., homotopy step).
Note that an eigenvalue crosses the unit circle around 5.6 ms as indicated by the gray dashed line. 
This eigenvalue actually crosses through $\lambda = -1$, and is indicative of a period doubling bifurcation.
The right hand plot in the figure shows a ``zoomed in'' version of Fig.~\ref{fig:circuit:bridge-failure}, which demonstrates that the period-doubling bifurcation precedes solver failure.
In particular, the time stepper successfully completes around 250 steps before oscillations are visible and another 250 before solver failure.

\begin{figure}[tp]
	\centering
	\includegraphics[width=0.45\textwidth]{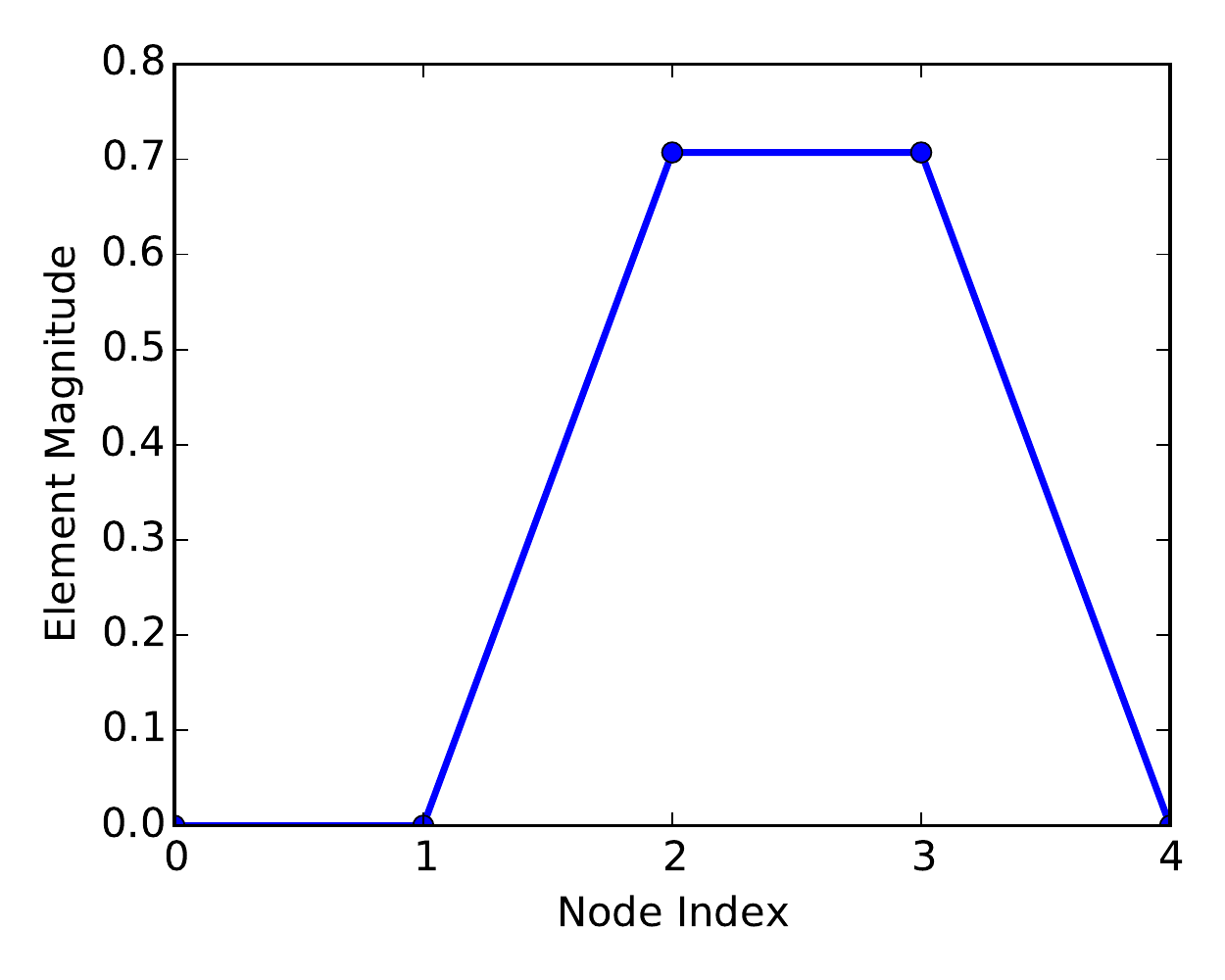}
	\includegraphics[width=0.45\textwidth]{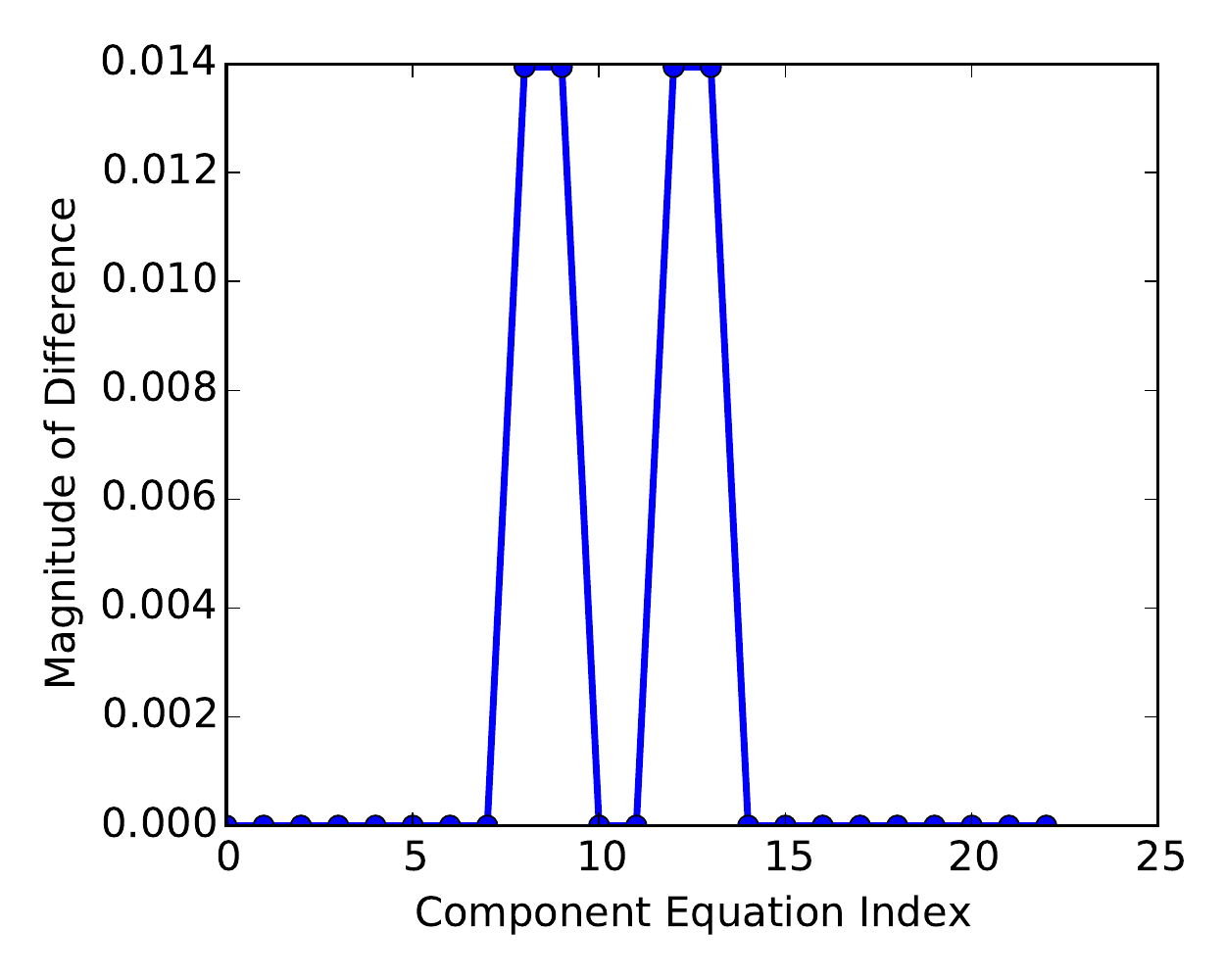}
	\caption{(left) The magnitude of the unstable eigenvector near the period-doubling bifurcation.
		(right) The difference between the implemented and finite-difference component-level directional derivatives.
		The indicated equations -- 8, 9, 12, and 13 -- are the component-level residual equations where implementation errors were introduced.
		There are two equations per error because there are two nodes associated with each of the two incorrectly implemented components.
	}
	\label{fig:circuit:bridge-failure-localization}
\end{figure}

As shown in this example, the eigenvalues indicate that ``something is going wrong'' with our the solver.
Except for the last step (which failed completely), none of the nonlinear solves required more than four Newton iterates, and most required only a pair of iterations.
At the time of failure, the system Jacobian is poorly conditioned, but that is also true in the reference example where issues with the nonlinear solver do not arise.

As shown in Fig.~\ref{fig:circuit:bridge-failure-localization}, the eigenvector indicates that nodes 2 and 3 are primarily affected by the loss of stability.
This corresponds to the time-trace shown in Fig.~\ref{fig:circuit:bridge-failure-comparison}, where oscillations in those nodes are visible just prior to solver failure.
To identify which components may be at fault, we compute the directional derivative of the residual function along this eigenvector.
A comparison between the implemented Jacobian and the finite-difference approximation reveals the convergence anomaly occurs in residual equations 8, 9, 12, and 13, which are the four equations with the Jacobian errors.
Therefore, the method was able to correctly identify and trace the error introduced into the system.

\subsubsection{A Single Component with an Implementation Error}

In the previous example, the presence of an error in the Jacobian eventually led to the total failure of the numerical solver.
However, this is not always the case.
In this example, we take the diode bridge shown in Fig.~\ref{fig:circuit:bridge-error}, but ``correct'' the error in the blue diode. 
As a result, there is a single component (the red diode) with a sign flip in the component-level Jacobian.

\begin{figure}[tp!]
	\centering
	\includegraphics[width=0.8\textwidth]{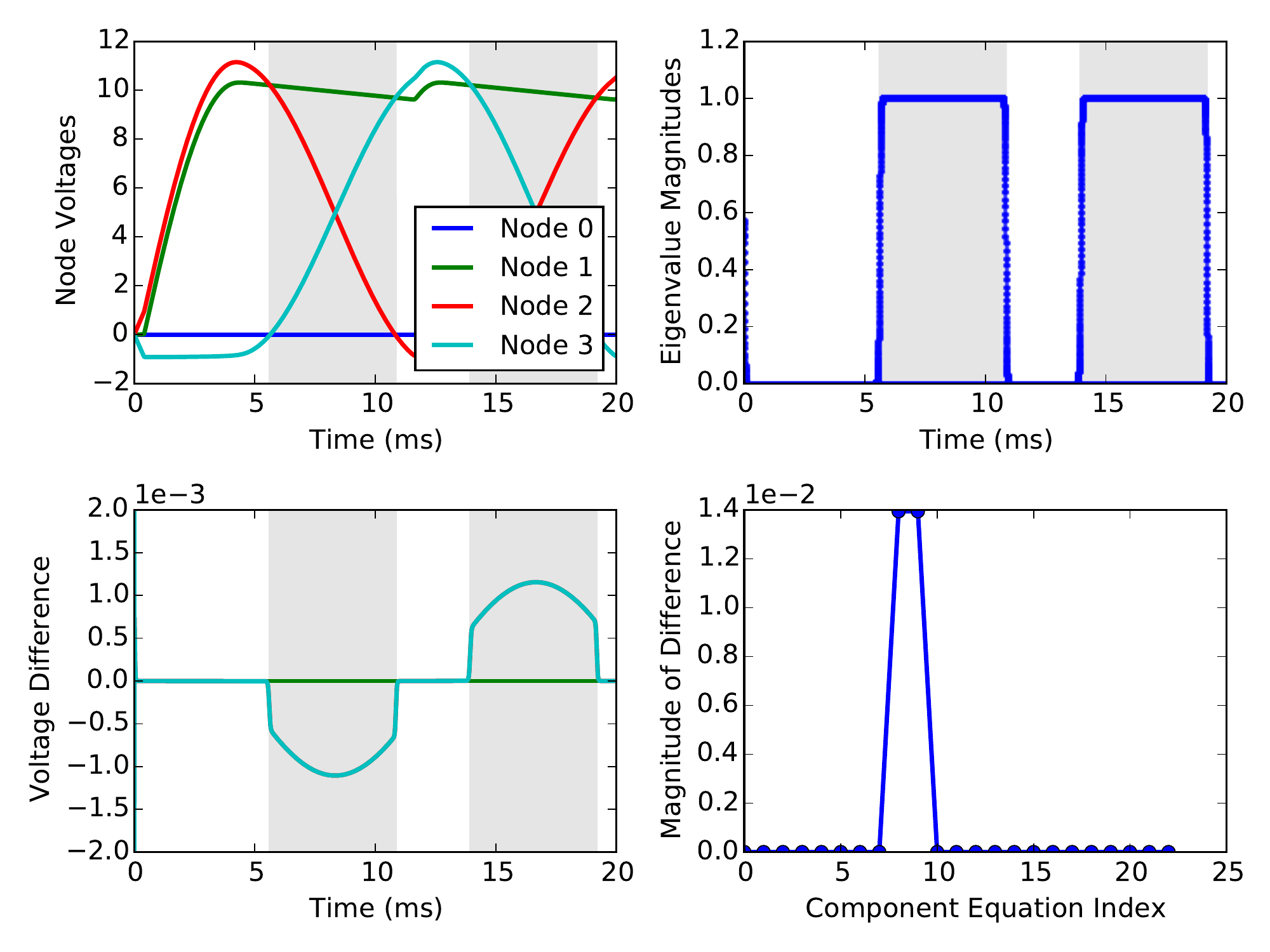}
	\caption{(top left) Plot of voltage versus time for the diode bridge with an implementation error in the diode indicated in red in Fig.~\ref{fig:circuit:bridge-error}.
		Despite the error, the nonlinear solver converges for every time step in the 20 ms window.
		(top right) Plot of the eigenvalue magnitudes as a function of time.
		(bottom left) The differences in the node voltages between the current and reference implementations; nodes 2 and 3 differ from their true values periodically in this time window.
		(bottom right) The eigenvector-based directional derivative check shows a difference between the implemented and true function values in component residual equations 8 and 9.
		These correspond to the equations where the sign flip in the component-level Jacobian was introduced.
		The shaded gray regions in the plots are the time intervals where at least one eigenvalue is greater than 0.5 in magnitude.	
	}
	\label{fig:circuit:bridge-converges}
\end{figure}

Figure~\ref{fig:circuit:bridge-converges} shows the trajectory and eigenvalue magnitudes associated with this system.
The results shown in this figure compare favorably with those in Fig.~\ref{fig:circuit:bridge-fine}.
Although not shown, the Newton solver converges after only a few iterations, and at least to the eye, the time-integration code appears to be functioning properly.
As shown in right-hand plot, there is a single eigenvalue that periodically ventures out to the unit circle and remains there for the intervals in time with the gray background.

The presence of an eigenvalue near the unit circle in this example does not indicate the solver will fail, but it does mark regions where our confidence in the solution should be diminished.
The bottom row of Figure~\ref{fig:circuit:bridge-converges} shows the difference between these trajectories and the reference implementation in  Fig.~\ref{fig:circuit:bridge-fine} as a function of time.
Note that a fixed time step integrator is being used, and that the nonlinear solver is Newton's method whose stopping criterion is when the norm of the residual is less than $10^{-8}$.
In short, the only difference between these two examples is the error introduced in the Jacobian.
The result, however, is a small but systematic discrepancy between the ``true'' and ``flawed'' solutions that appear in the gray shaded regions where a large eigenvalue exists.
In this example, the discrepancy vanishes outside of these regions, but in another application this error could grow as the simulation progresses.
The directional derivative check shown in the bottom-right plot localizes the error to equations 8 and 9, which are the two equations where the errors were introduced.

\begin{figure}[tp!]
	\centering
	\includegraphics[width=0.75\textwidth]{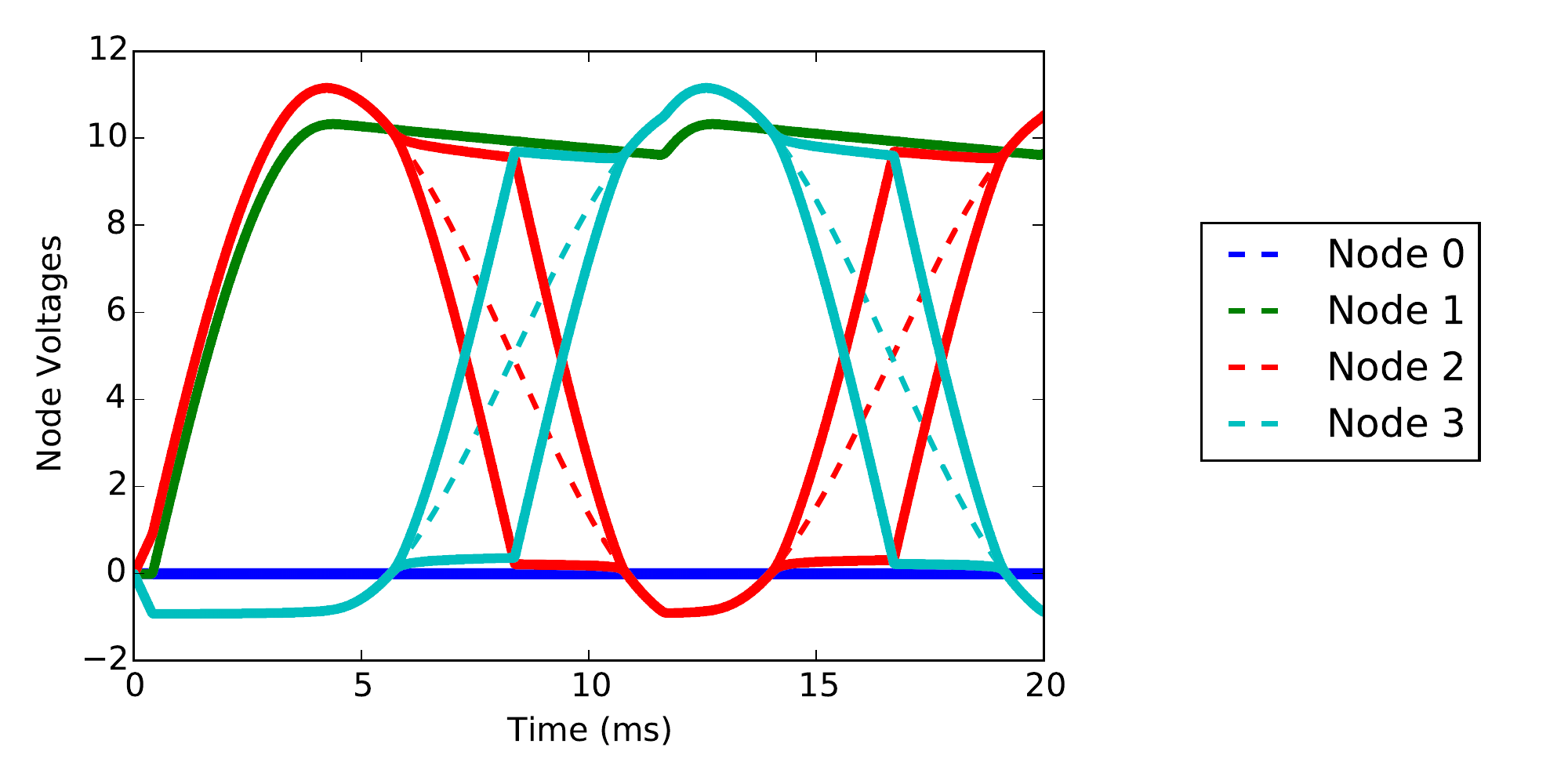}
	\caption{Plot of the Newton iterates at every time step when the nonlinear solver is forced to run twenty for iterations per step.
		For the sake of comparison, the dashed lines reproduce the trajectories obtained when a residual norm of $10^{-8}$ is used as the stopping criterion.
		The predictions of the two methods diverge in these regions, with the period doubling bifurcation causing the values of nodes 2 and 3 to jump between two possible sets of values.
		The small but systematic errors highlighted in Fig.~\ref{fig:circuit:bridge-converges} are due to the fact the nonlinear solver is ``getting lucky'' and finding a solution that meets the residual tolerance rather than converging to a fixed point.
	}
	\label{fig:circuit:bridge-converges-period-doubling}
\end{figure}

The cause of the discrepancy between the true solution and the one introduced here are due to the interplay of the solver dynamics and the stopping criterion.
At the start and end of the regions of interest, a period-doubling bifurcation occurs.
As a result, Newton's method will no longer converge to a single point, but will instead oscillate between a pair of points.
The transition from a single solution to a period-two orbit is shown in Fig.~\ref{fig:circuit:bridge-converges-period-doubling} where the stopping criterion was removed.
The solid lines show the iterations taken by Newton's method at each time step, and the dashed lines show the trajectories obtained in Fig.~\ref{fig:circuit:bridge-converges}.
Note how the nonlinear solver switches between identifying a single solution to alternating between two choices at different intervals in time.

The loss of stability of the original branch has a profound impact on the confidence one should have in the resulting solution.
As shown here, the nonlinear solver is not converging to a single fixed point, and the results we are observing in Fig.~\ref{fig:circuit:bridge-converges} are due to transient points that have residuals whose norm is below the tolerance level of $10^{-8}$.
Therefore, seemingly minor changes to the code, such as lowering the tolerance from $10^{-8}$ to $10^{-10}$, which can be done without issue in the original system, or simply requiring a minimum number of Newton iterations to be performed can result in solver failure.
As a result, even though the time-integration procedure did not fail and the resulting solutions resemble the correct ones, the presence of a bifurcation signals that one should lose confidence in the resulting solution because it is unclear what of the many possible behaviors (e.g., switching to one of many branches, stopping at a transient point, etc.) is responsible for the resulting solutions.

\subsubsection{Example: Parameter Dependent Instabilities in a Diode Bridge}

Finally, we study the effects of removing the large parallel resistors that implement the minimum conductances of the diodes. 
In the previous section, the value of this resistor was $10^{12}$ $\Omega$, and for this example we will steadily increase it up to $10^{35}$ $\Omega$. 
In principle, as the resistance approaches infinity, the parallel resistor is effectively removed and the diode acts ideally, but as we will demonstrate shortly, our time-integration routine will exhibit numerical issues at large but finite values of parallel resistance.

\begin{figure}[tp!]
	\centering
	\includegraphics[width=0.7\textwidth]{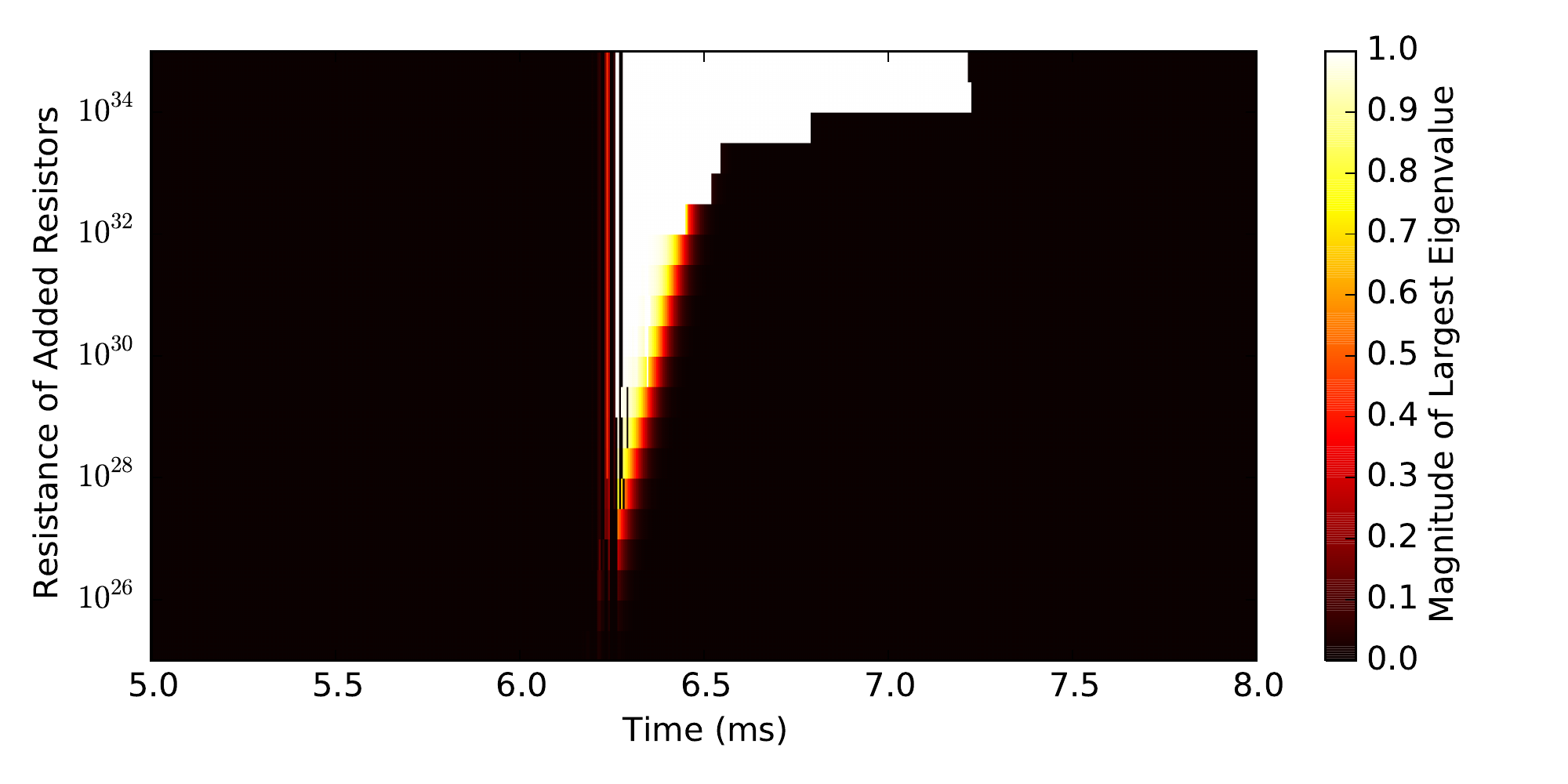}
	\caption{Plot of the magnitude of the largest eigenvalue as the minimum conductances, which are equal to the reciprocal of the parallel resistance, is removed.
		Note that as the parallel resistance grows, the time-interval where the eigenvalue is large also grows.
		Despite this eigenvalue, the solver converges for all time steps points on the plot above. 
	}
	\label{fig:circuit:bridge-singularity}
\end{figure}

Figure~\ref{fig:circuit:bridge-singularity} shows the magnitude of the largest eigenvalue at various time steps and values of the parallel resistor.
To generate this figure, multiple experiments were performed with fixed values of the parallel resistance, so the pseudo-color plot is comprised of multiple ``sweeps'' in time at fixed values of resistance.
The main feature of interest is the white region, which corresponds to combinations of time and resistance with an eigenvalue on or near the unit circle. 
As in the previous example, this does not result in solver failure, but it does have an impact on the resulting solution.

\begin{figure}[tp!]
	\centering
	\includegraphics[width=0.45\textwidth]{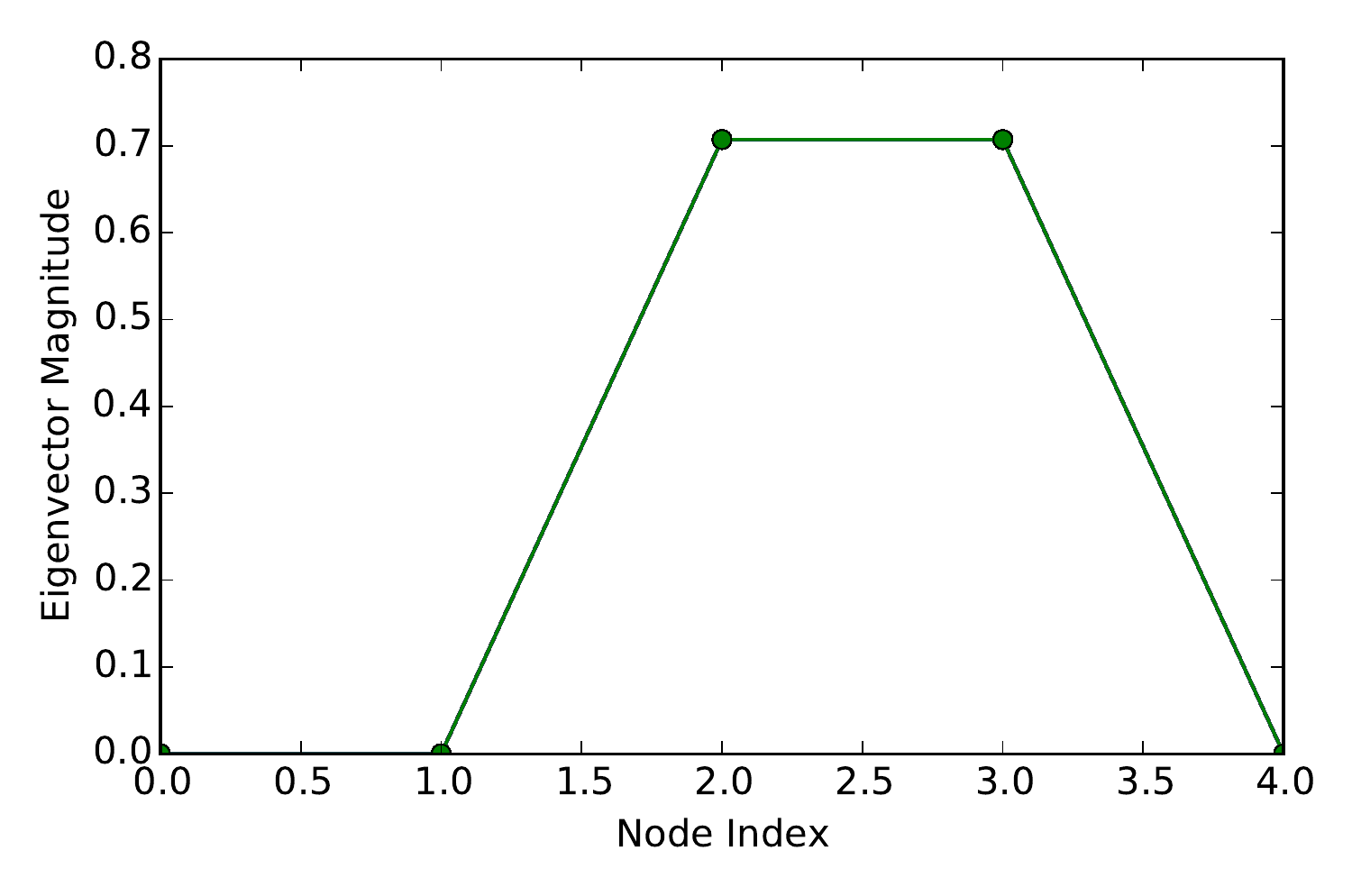}
	\caption{Plot of the eigenvectors associated with an eigenvalue near the period doubling bifurcation for the entire range of added resistances shown in Fig.~\ref{fig:circuit:bridge-singularity}.
		The corresponding eigenvector highlights that nodes 2 and 3 are affected by the increase in the parallel resistances, while the other three nodes are not.
	}
	\label{fig:circuit:bridge-singularity-localization}
\end{figure}

Figure~\ref{fig:circuit:bridge-singularity-localization} shows the superposition of the eigenvectors associated with the least-stable eigenvalue at the first and last points of the ``white region'' in Fig.~\ref{fig:circuit:bridge-singularity} over all of the resistances. 
As shown in this figure, only the second and third nodes are affected by the convergence  anomaly that appears as the parallel resistors are removed (i.e., as $R\to\infty$).
The other three nodes appear to be unaffected by this change, as they do not appear in any significant way in the resulting eigenvector.
Applying the directional derivative check used in Sec~\ref{sec:theory:localize} does not result in any significant mismatches between the true and numerically implemented system- and component-level Jacobians.
% TODO: change system-wide to something else
As a result, the stability analysis for this particular example indicated a system-wide problem rather than a component-specific issue.

\begin{figure}[tp!]
	\centering
	\includegraphics[width=\textwidth]{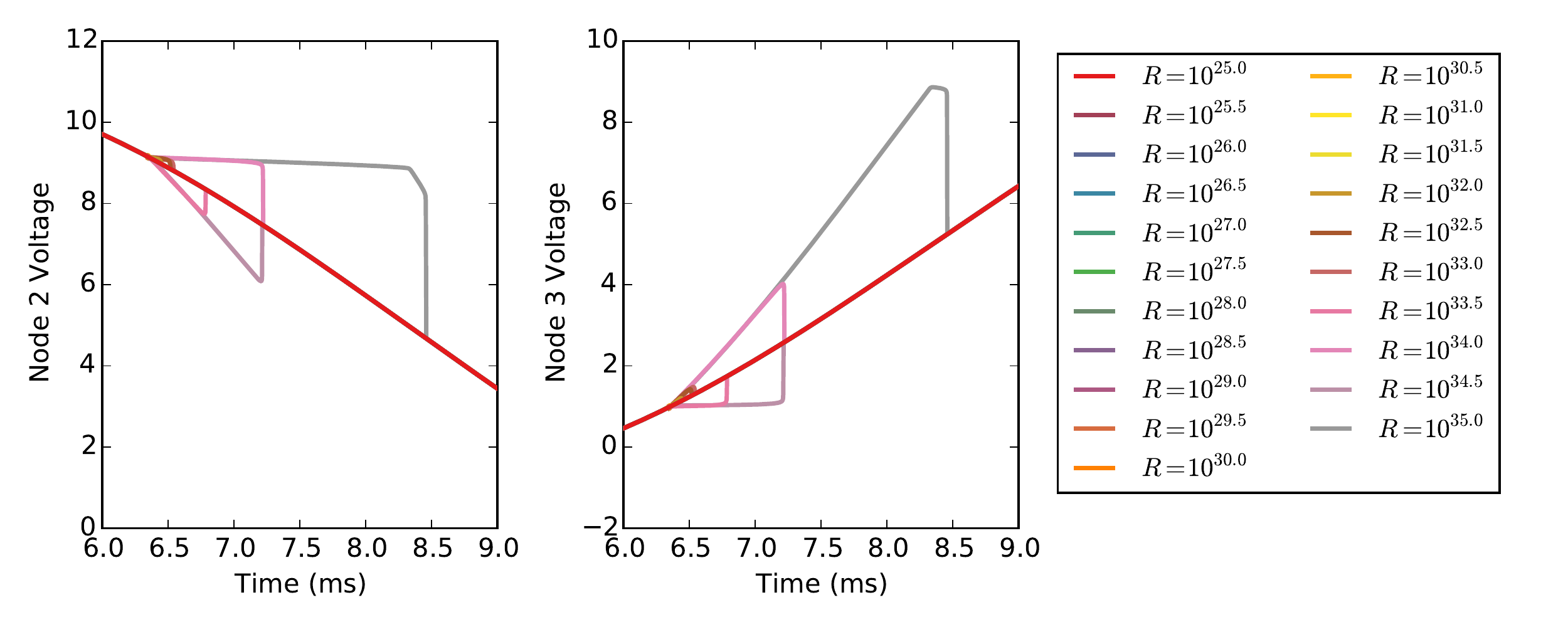}
	\caption{Plot of the voltages in nodes 2 and 3 over the range of parallel resistors shown in Fig.~\ref{fig:circuit:bridge-singularity}.
		In the interval where an eigenvalue is near one in magnitude, these voltages are free to ``float,'' which results in the large differences in voltage observed in these nodes in these regions.
	}
	\label{fig:circuit:bridge-singularity-float}
\end{figure}

To highlight the impact of the bifurcation, Figure~\ref{fig:circuit:bridge-singularity-float} shows the voltage trajectories of nodes 2 and 3, which are the two nodes affected by the removal of the parallel resistors. 
The other three nodes are not impacted in any meaningful way, and their evolution can be seen in Fig.~\ref{fig:circuit:bridge-fine}.
Note the voltages in nodes 2 and 3 have large ``excursions'' away from the trajectories that result with smaller resistances, which persist as long as the eigenvalue remains near the unit circle.
Intuitively, these systems should not be physically different -- all cases have a very large resistor in parallel with a diode -- so this difference is numerical in nature.

The behavior shown in Fig.~\ref{fig:circuit:bridge-singularity-float} is similar to the behavior shown in Fig.~\ref{fig:circuit:bridge-converges-period-doubling}.
As the parallel resistors are removed, the original solution branch loses stability and Newton's method begins to oscillate between two points.
However, the residuals associated with these points can be small enough to meet the convergence criterion, and a point is accepted as the next time step.

Physically these results can be explained.
When the parallel resistors are too large, nodes 2 and 3 effectively lose connection to the ground.
As this is an electrical system, the only quantity of physical importance is the {\em differences} in voltages, so without a connection to the ground, nodes 2 and 3 can take on any value they please provided their difference matches the constraints imposed by the voltage source.
However, the value of resistance that is ``too large'' depends on time, so the interval where nodes 2 and 3 are ``floating'' are affected by the parallel resistors.

\subsection{Example: Multiple Simultaneous Errors in a Power Channel}
\label{[sec:power-channel]}

\begin{figure}[tp!]
	\centering 
	
	\usetikzlibrary{circuits.logic.US,circuits.logic.IEC}
	\scalebox{0.45}{
		\begin{tikzpicture}[every node/.style={outer sep=0, inner sep=0},  line width=0.3mm]
		
		\node[inner sep=0] (n0) at (0,0) {};
		\node (n1) at (0,3) {};
		
		\draw (n0)  to[sI=\parbox{0.75in}{\centering\huge 12 V \\ 60 Hz}] (n1);
		
		\node (v1) at (3.5,3) {};
		\node(n2) at (5.5,5) {};
		\node(n3) at (5.5,1) {};
		\node (v2) at (7.5,3) {};
		
		\draw (n3) to[diode, *-*] (v2);
		\draw (v1) to[diode, color=red] (n2);
		\draw[ultra thick, color=red] ($0.5*(v1)+0.5*(n2)$) circle (0.6cm);
		\draw (v2) to[diode,  *-*] (n2);
		\draw (n3) to[diode,  *-*] (v1);
		\node (v3) at (9,5) {};
		\node at (9,1) {};
		\node (v4) at (11,5) {};
		\node at (11,1) {};
		\draw (n2) to[short,  *-*] (v3);
		\draw (v3) to[short,  *-*] (v4);
		\draw (11,5) to[resistor={\hspace{.1in}\parbox{0.75in}{\centering\huge 20 $\Omega$ \\ 5 mF}},  *-*] (11,1) node (v6) {};
		\draw (v3) to[capacitor,  *-*] (9,1) node (v5) {};
		\draw (n3) to[short,  *-*] (v5);
		\draw (v5) to[short,  *-*] (v6);
		
		\node (v7) at (2,3) {};
		\draw (n1) to[short,  *-*] (v7);
		\draw (v7) to[short,  *-*] (v1);
		
		\node (vb1) at (3.5,-2) {};
		\node(nb2) at (5.5,0) {};
		\node(nb3) at (5.5,-4) {};
		\node (vb2) at (7.5,-2) {};
		
		\draw (nb3) to[diode,  *-*] (vb2);
		\draw (vb1) to[diode,  *-*] (nb2);
		\draw (vb2) to[diode,  *-*] (nb2);
		\draw (nb3) to[diode,  *-*] (vb1);
		\node (vb3) at (9,0) {};
		
		\node (vb4) at (11,0) {};
		\draw (nb2) to[short,  *-*] (vb3);
		\draw (vb3) to[short, *-*] (vb4);
		\draw (11,0) to[resistor={\hspace{.1in}\parbox{0.75in}{\centering\huge 2 k$\Omega$ \\ 5 mF}},  *-*] (11,-4) node (vb6) {};
		\draw (vb3) to[capacitor,  *-*] (9,-4) node (vb5) {};
		\draw (nb3) to[short,  *-*] (vb5);
		\draw (vb5) to[short,  *-*] (vb6);
		
		\node (vc1) at (3.5,-7) {};
		\node(nc2) at (5.5,-5) {};
		\node(nc3) at (5.5,-9) {};
		\node (vc2) at (7.5,-7) {};
		
		\draw (nc3) to[diode,  *-*] (vc2);
		\draw (vc1) to[diode, *-*] (nc2);
		\draw (vc2) to[diode, *-*] (nc2);
		\draw (nc3) to[diode, *-*] (vc1);
		\node (vc3) at (9,-5) {};
		
		\node (vc4) at (11,-5) {};
		\draw (nc2) -- (vc3);
		\draw (vc3) -- (vc4);
		\draw (11,-5) to[resistor={\hspace{.1in}\parbox{1in}{\centering\huge 10 m$\Omega$ \\ 5 mF}}, *-*] (11,-9) node (vc6) {};
		\draw (vc3) to[capacitor, *-*] (9,-9) node (vc5) {};
		\draw (nc3) to[short,  *-*] (vc5);
		\draw (vc5) to[short,  *-*] (vc6);
		
		\node (vd1) at (3.5,-12) {};
		\node(nc2) at (5.5,-10) {};
		\node(nc3) at (5.5,-14) {};
		\node (vd2) at (7.5,-12) {};
		
		\draw (nc3) to[diode, *-*] (vd2);
		\draw (vd1) to[diode, *-*] (nc2);
		\draw (vd2) to[diode, *-*] (nc2);
		\draw (nc3) to[diode, *-*] (vd1);
		\node (vd3) at (9,-10) {};
		
		\node (vd4) at (11,-10) {};
		\draw (nc2) -- (vd3);
		\draw (vd3) to[short, *-*] (vd4);
		\draw (11,-10) to[resistor={\hspace{.1in}\parbox{0.75in}{\centering\huge \vspace*{3em} 10 $\Omega$ \\ 10 $\mu$F \\ 1 mH}}, *-*] (11,-12);
		\draw (11,-12) to[inductor,color=blue] (11,-14) node (vd6) {};
		\draw[ultra thick, densely dashed, color=blue] (11, -13) circle (0.6cm);
		\draw (vd3) to[capacitor] (9,-14) node (vd5) {};
		\draw (nc3) -- (vd5);
		\draw (vd5) to[short,  *-*] (vd6);
		
		\node (v13) at (7.5,6) {};
		\node (v12) at (2,10) {};
		\node (v8) at (9,10) {};
		\node (v9) at (9,6) {};
		\node (v10) at (11,10) {};
		\node (v11) at (11,6) {};
		\draw (v8) to[capacitor, *-*] (v9);
		\draw (v10) to[resistor={\hspace{.1in}\parbox{0.75in}{\centering\huge 1 $\Omega$ \\ 5 mF}}, *-*] (v11);
		\draw (v10) -- (v8) to[short, *-*] (v12);
		\draw (v11) -- (v9) -- (v13);
		
		\node (v19) at (7.5,11) {};
		\node (v14) at (2,15) {};
		\node (v15) at (8,15) {};
		\node (v16) at (8,11) {};
		\node (v17) at (9.5,15) {};
		\node (v18) at (9.5,11) {};

		\draw (v14) to[diode, *-*] (v15);
		\draw (v15) to[capacitor, *-*] (v16);
		\draw (v17) to[resistor, *-*] (v18);
		\draw (v15) -- (v17);
		\draw (v16) -- (v18);
		\draw (v19) -- (8,11);

		\node (ve14) at (9.5,15) {};
		\node (ve15) at (12,15) {};
		\node (ve16) at (12,11) {};
		\node (ve17) at (13.5,15) {};
		\node (ve18) at (13.5,11) {};
		
		\draw (ve14) to[diode, *-*] (ve15);
		\draw (ve15) to[capacitor, *-*] (ve16);
		\draw (ve17) to[resistor, *-*] (ve18);
		\draw (ve15) -- (ve17);
		\draw (ve16) -- (ve18);
		
		\node (vf19) at (17.5,11) {};
		\node (vf14) at (13.5,15) {};
		\node (vf15) at (16,15) {};
		\node (vf16) at (16,11) {};
		\node (vf17) at (17.5,15) {};
		\node (vf18) at (17.5,11) {};
		
		\draw (vf14) to[diode, *-*] (vf15);
		\draw (vf15) to[capacitor, *-*] (vf16);
		\draw (vf17) to[resistor={\hspace{.1in}\parbox{0.75in}{\centering\huge 1 $\Omega$ \\ 1 pF}}, *-*] (vf18);
		\draw (vf15) -- (vf17);
		\draw (vf16) -- (vf18);
		\draw (vf19) -- (8,11);

		\node (v20) at (2,-2) {};
		\draw (v20) to[short, *-*] (vb1);
		\node (v21) at (2,-7) {};
		
		\node[inner sep=0] (v22) at (2,-12) {};
		
		\draw (v21) to[short, *-*] (vc1);
		\draw (v22) -- (vd1);
		\draw (v14) -- (2,10) -- (v7) -- (2,-2) -- (v21) -- (v22);
		\draw[color=gray] (v19) -- (v13) -- (v2) -- (vb2) -- (vc2) -- (vd2) -- (7.5,-15) -- (0,-15) node (v23) {};
		\draw[color=gray] (n0) -- (v23);
		\draw (v23) node[ground, color=gray]{};
		\end{tikzpicture}
	}
	\caption{Circuit diagram of a simple power channel. 
		A single AC source drives six loads (3 DC loads and 3 AC loads).
		The size of the components in each load are given on the right; for example, every resistor in the top-most load is 1 $\Omega$ and every capacitor is 1 pF.
		The diodes are implemented with a parallel resistor of $10^{12}$ $\Omega$, which is not shown in the diagram above. 
		The resistor in red (solid circle) and the inductor in blue (dashed circle) have implementation errors in their Jacobians.
		In particular, the diode in red has a sign flip in its Jacobian, and the inductor has a Jacobian where one element is multiplied by a factor of 0.95.
	}
	\label{fig:circuit:power-channel}    
\end{figure}
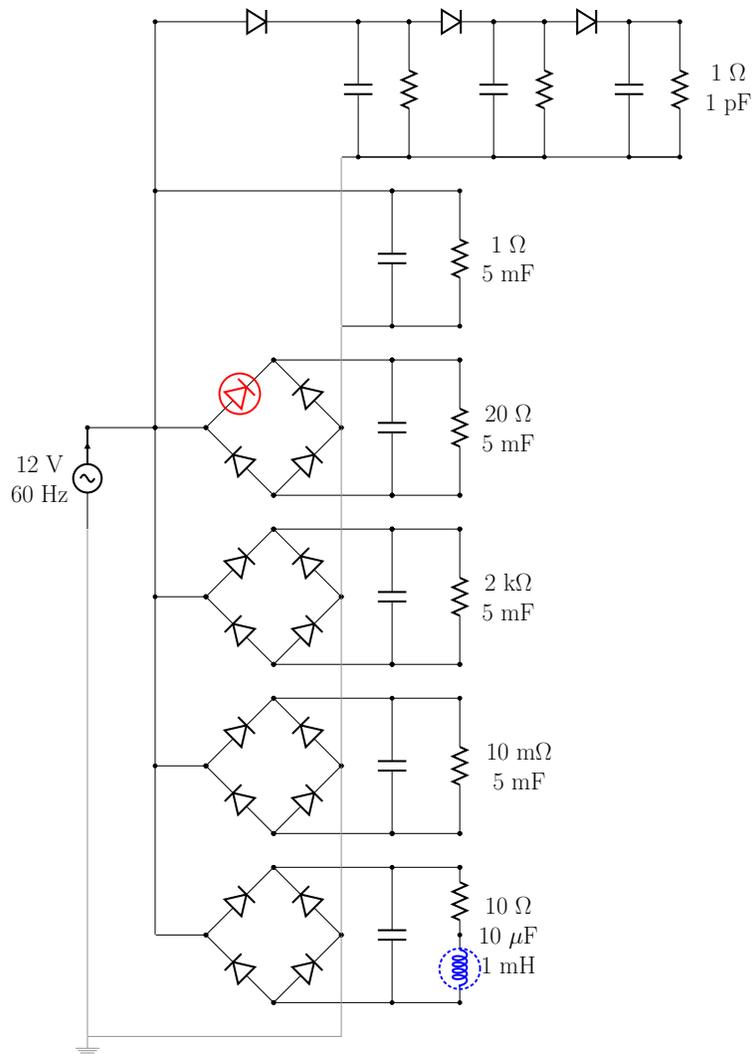

The final circuit-based example is an expanded version of the diode bridge, and has the schematic shown in Fig.~\ref{fig:circuit:power-channel}.
This example consists of an ideal AC source driving six different loads.
From top to bottom they are:
\begin{enumerate}
	\item a diode chain (see, e.g,~\cite{voss2008model});
	\item a resistor and a capacitor in parallel;
	\item a resistive load of 20 $\Omega$ behind a diode bridge;
	\item a relatively large resistive load of 2 k$\Omega$ behind a diode bridge;
	\item a relatively small resistive load of 10 m$\Omega$ behind a diode bridge;
	\item a resistor and nonlinear inductor behind a diode bridge; note that the filter capacitor is undersized, and therefore does not smooth out the ripples in the bridge output.
\end{enumerate}
This circuit contains two types of nonlinear components: diodes and a nonlinear inductor whose residual equation is:
\begin{equation}
R_L(V_0, V_1, i_L) = 
\begin{bmatrix}
i_L \\ -i_L \\ V_0 - V_1 - \frac{\partial\lambda}{\partial i_L}\frac{di_L}{dt}
\end{bmatrix},
 \qquad \lambda(i_L) = \frac{L_0I_{\text{sat}}}{\sqrt{I_{\text{sat}}^2 + i_L^2}}i_L,
\end{equation}
where $V_0$ and $V_1$ are the voltages on either terminal of the inductor, and $i_L$ is an internal state containing the current through the nonlinear inductor.
The parameters $L_0 = 0.001$ H and $I_{\text{sat}} = 1$ A are the ``linear'' inductance and saturation current of the nonlinear inductor.
We will use this system to benchmark our approaches on two types of issues: multiple implementation errors and step size selection for the BDF.
The first is a component-level issue, which the approach can isolate, while the second is systemic and appears when the step size becomes too small.

\subsubsection{Diagnosing Multiple Implementation Errors}

There are two implementation errors that will be introduced to this system, and appear in the components shown in red and blue (solid and dashed circles respectively).
The diode has a sign flip in the component-level Jacobian, and is the same error we identified in the single diode bridge.
The inductor possesses a term that is scaled by a factor of 0.95, which was implemented by using two different values of the inductance when computing the residual and the component-level Jacobian.
Our objective in this example is to demonstrate that the localization procedure still works in a larger system when: (a) multiple errors are present, and (b) some sub-systems are error free. 

\begin{figure}[tp!]
	\centering
	\includegraphics[width=0.8\textwidth]{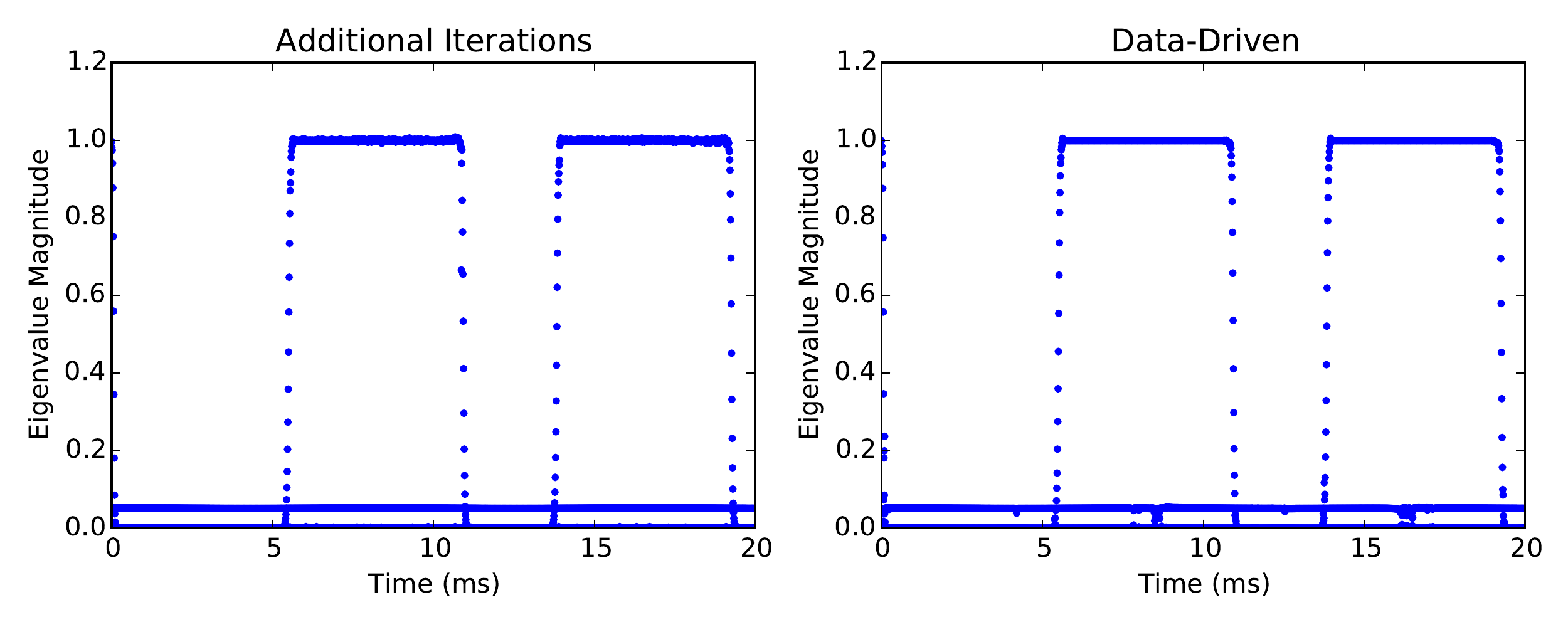}
	\caption{Plot of the eigenvalues of the linearized Newton-solver at each time step of the system shown in Fig.~\ref{fig:circuit:power-channel}.
		The left plot shows the eigenvalues obtained using the additional iterations approach, and the right plot shows the eigenvalues obtained using the data-driven approach.
		While the eigenvalues and eigenvectors obtained from the data-driven approach will be used in the analysis that follows, both methods identify a persistent ``smaller'' eigenvalue with $|\lambda|\approx 0.05$ and an intermittent larger eigenvalue of $\lambda \approx -1$.
	}
	\label{fig:circuit:power-channel-eigenvalues}
\end{figure}

We integrated this system over a 20 ms time window using backward Euler with a time step of $\Delta t = 10^{-3}$ ms.
A standard Newton solver was used to solve the nonlinear system of equations that results at each time step.
The initial guess was the solution at the previous time step, and the convergence criterion was the norm of the residuals dropping below $10^{-8}$.
As this is a large system with convergence anomalies, a larger number of Newton iterates were required at each time step, and the number of iterations typically ranged between 5 and 10.

Because we have more data, the ``data-driven'' approach for approximating the eigenvalues of the linearized Newton solver produces useful results at every time step.
Figure~\ref{fig:circuit:power-channel-eigenvalues} shows the eigenvalues obtained using the ``additional experiments'' and ``data-driven'' approaches at each time step.
Note that there are now {\em two} convergence anomalies: one that remains constant with $|\lambda|\approx 0.05$ and the other that appears periodically and ventures out to $|\lambda|=1$.
There is good agreement between the two methods, so in the results that follow, the eigenvectors and eigenvalues will be those computed using the data-driven approach, though either method will produce similar results.

\begin{figure}[tp]
	\centering
	\includegraphics[width=0.45\textwidth]{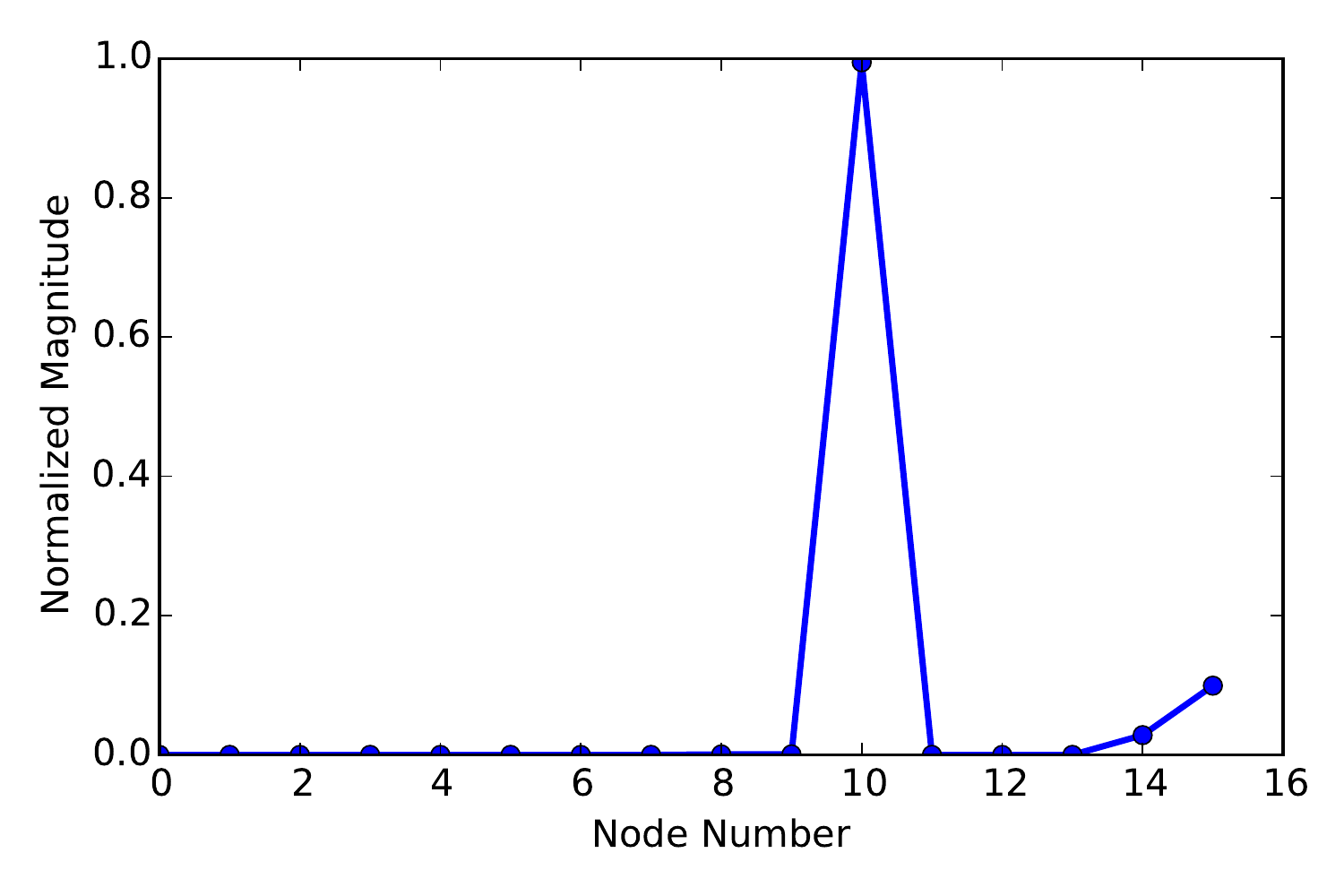}
	\includegraphics[width=0.45\textwidth]{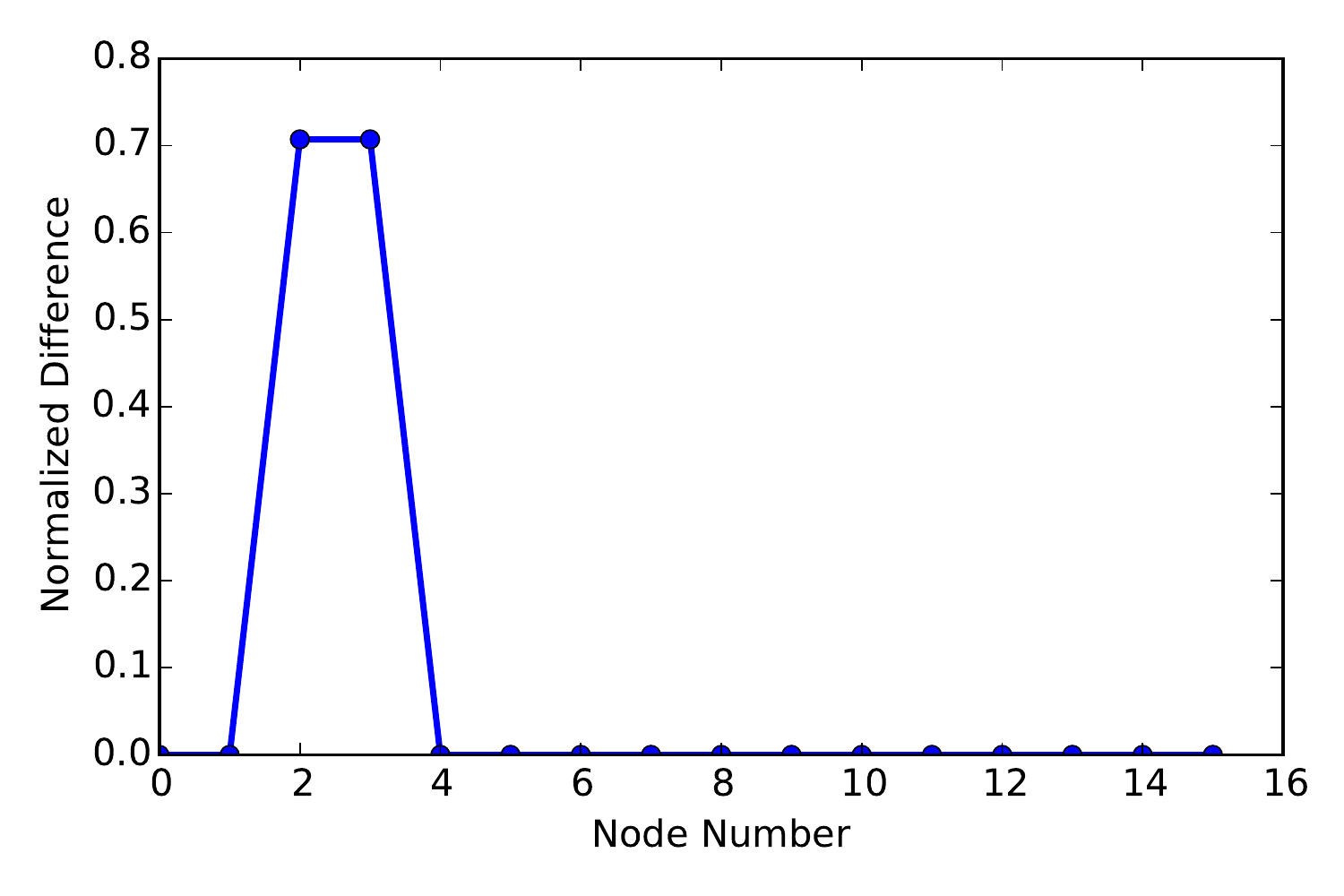}
	
	\includegraphics[width=0.45\textwidth]{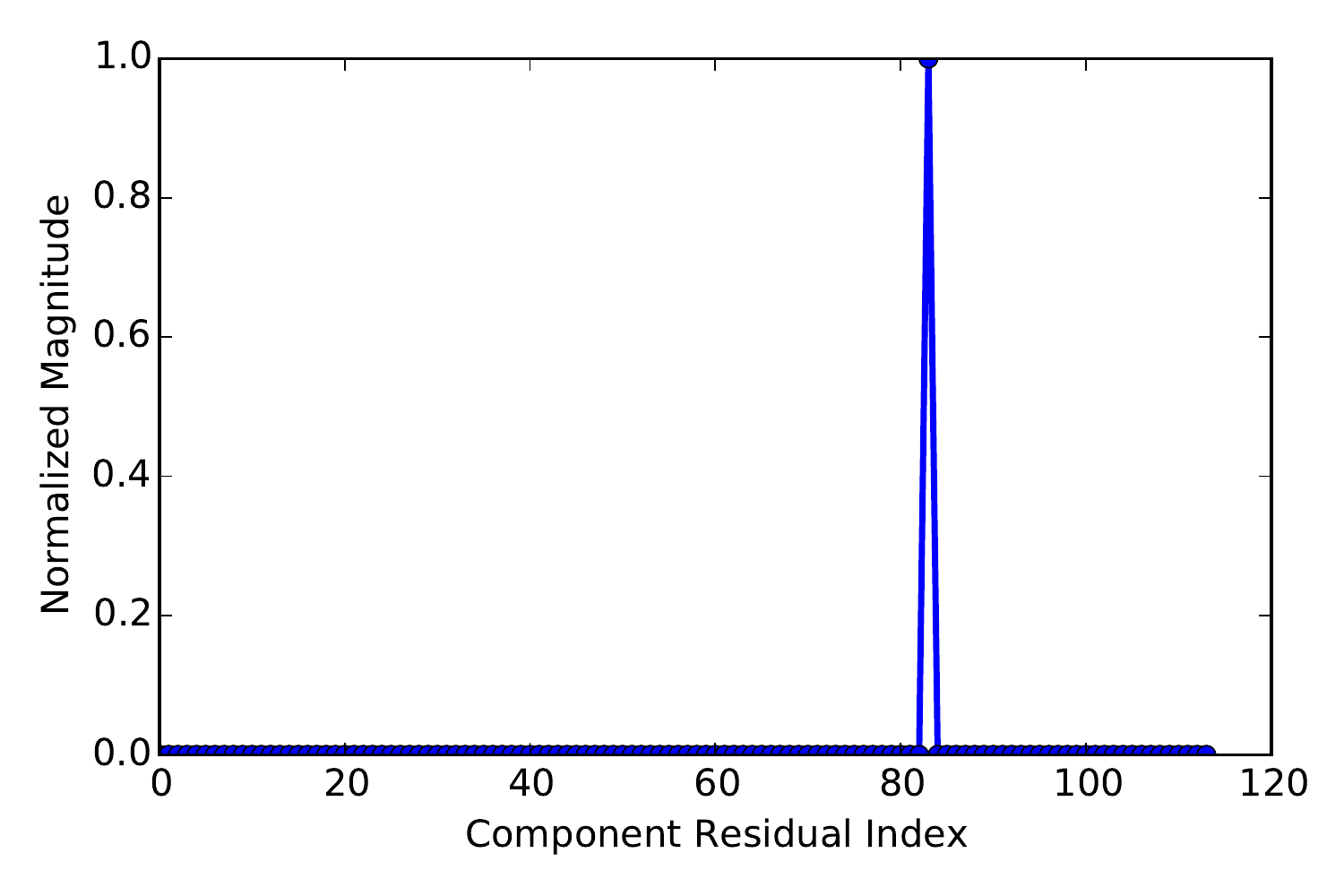}
	\includegraphics[width=0.45\textwidth]{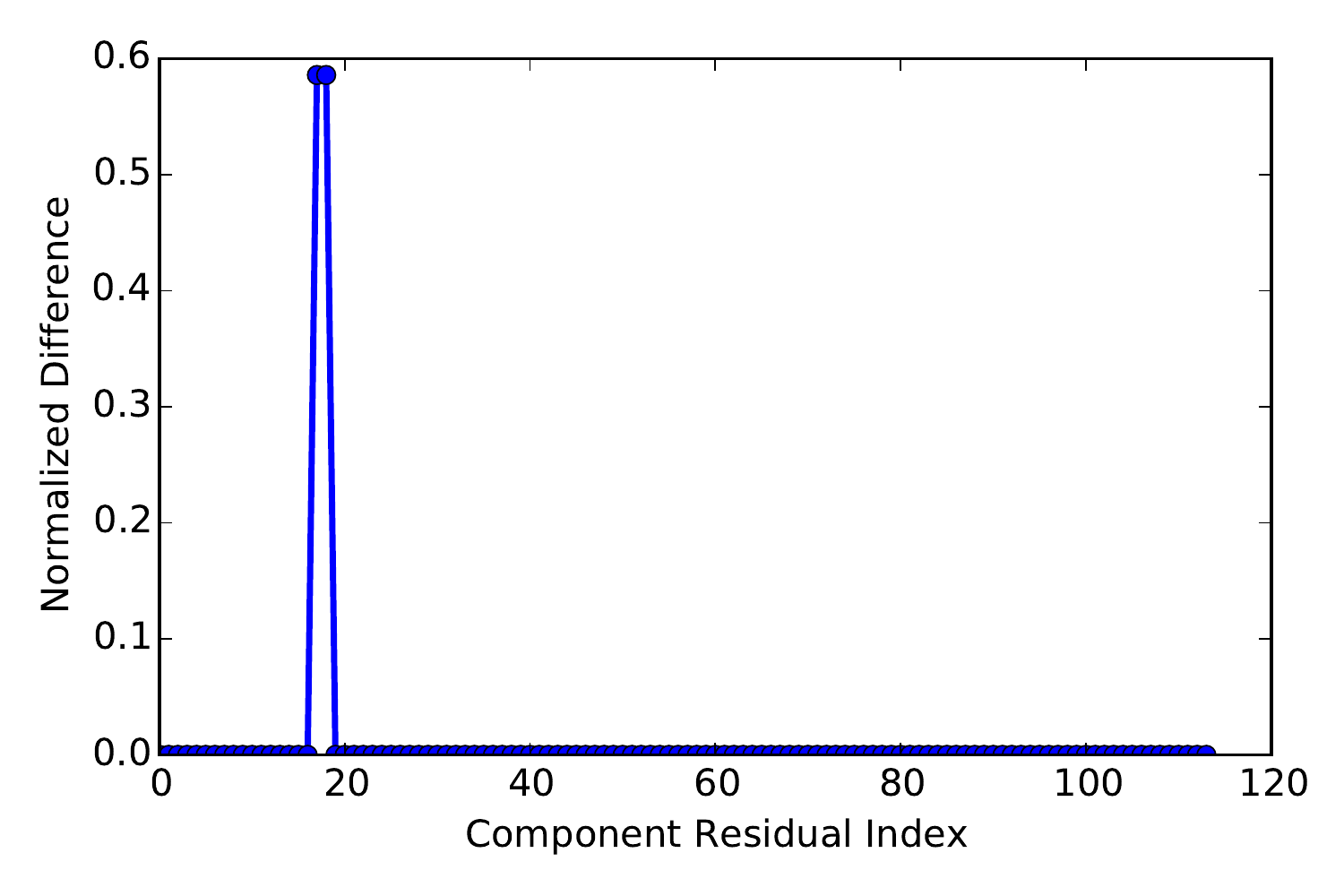}
	\caption{(top left) The eigenvector associated with $\lambda \approx 0.05$ obtained using the data-driven approach near $t=20$ ms. 
		(top right) The eigenvector associated with $\lambda \approx -1$ obtained using the data-driven approach near $t=19$ ms.
		Both eigenvectors identify the {\em nodes} that will converge to the fixed point slowly due to the presence of convergence anomalies. 
		(bottom) The differences between the implemented and finite-difference approximation of the component-level directional derivatives along the eigenvectors shown above.
		This analysis highlights component-level residual 83, which is the nonlinear inductor, and component-level residuals 17 and 18, where a sign flip in the resistor Jacobian was introduced.
	}
	\label{fig:circuit:power-channel-vectors}
\end{figure}

Figure~\ref{fig:circuit:power-channel-vectors} (top) shows the eigenvectors associated with the two largest eigenvalues.
Note that the eigenvectors pick out different nodes in the system: the eigenvector associated with $\lambda\approx 0.05$ is supported mainly on node 10, which is an input to the inductor, and nodes 14 and 15, which are internal nodes associated with the nonlinear inductor and voltage source.
On the other hand, the eigenvector with $\lambda\approx -1$ identifies nodes 2 and 3, which are the two nodes in the subsystem with the flawed diode.
The eigenvectors only indicate which nodes converge slowly to their final values, and not necessarily indicate the inputs to a flawed component.
In particular, no implementation errors were added to the voltage source, but its internal node still appeared as part of the eigenvector.
Similarly, node 3 is not physically connected to the flawed diode, but is prominent in the eigenvector.
Both nodes appear in the eigenvector because they are affected by the error, and not because they are the source of the issue.

What the corresponding eigenvectors allow us to do is identify discrepancies between the true and implemented Jacobians, which will can pinpoint where implementation errors are present.
The (scaled) difference in the directional derivative is shown in Fig.~\ref{fig:circuit:power-channel-vectors}.
The directional derivative check associated with $\lambda = 0.05$ highlights component equation 83, which is the equation where the scaling factor was introduced in the Jacobian of the nonlinear inductor.
Similarly, comparing the directional derivatives in the direction of the eigenvector with $\lambda \approx -1$ identifies equations 17 and 18, which are the residuals associated with the diode whose Jacobian had a sign error.

In this section, we demonstrated the ``additional experiments'' and ``data-driven'' approaches on a set of problems involving the simulation of electrical circuits.
In these examples, we diagnosed deliberately introduced issues in component-level equations in order to show both approaches could accurately detect the issue, and localize it to the affected nodes and, in this setting, the responsible component-level residual equations.
We also showed that the approach could detect system wide issues such as the effective loss of solution uniqueness that occurs in these simulations when the resistors that implement the minimum conductances of the diodes approach infinity.
In this case, the eigenvalues of the system linearization detected the issue and the eigenvectors isolated the problem to the affected nodes, but no further localization was possible.

\subsubsection{Example: The Effects of Time-Step Size on Solver Performance}
\label{sec:step-size}

Now we consider the effects that the choice of time step, $\Delta t$, have on the performance of the nonlinear solver.
Many numerical codes, such as the MATLAB implementation of embedded Runge-Kutta methods in \texttt{ode45}~\cite{dormand1980family}, attempt to bound the integration error by adapting the step size.
In particular, if the estimated error exceeds a certain threshold, then the step size is decreased.
However for systems governed by DAEs or implicit time steppers, the choice of the time step can have a measurable impact on the performance of the resulting solver. 

\begin{figure}[tp]
	\centering
	\includegraphics[width=0.45\textwidth]{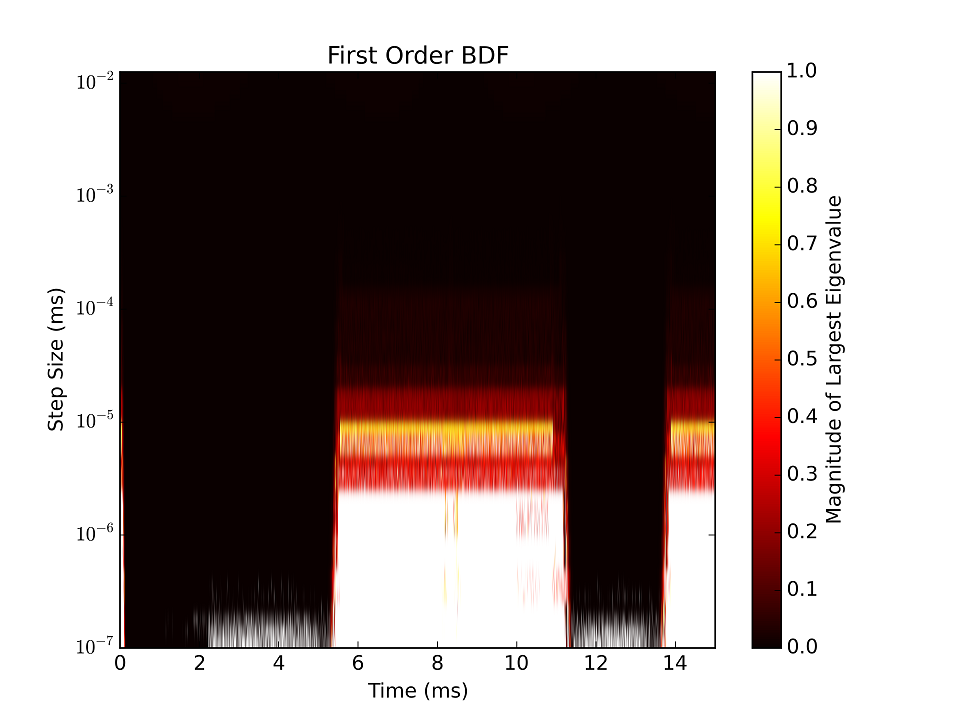}
	\includegraphics[width=0.45\textwidth]{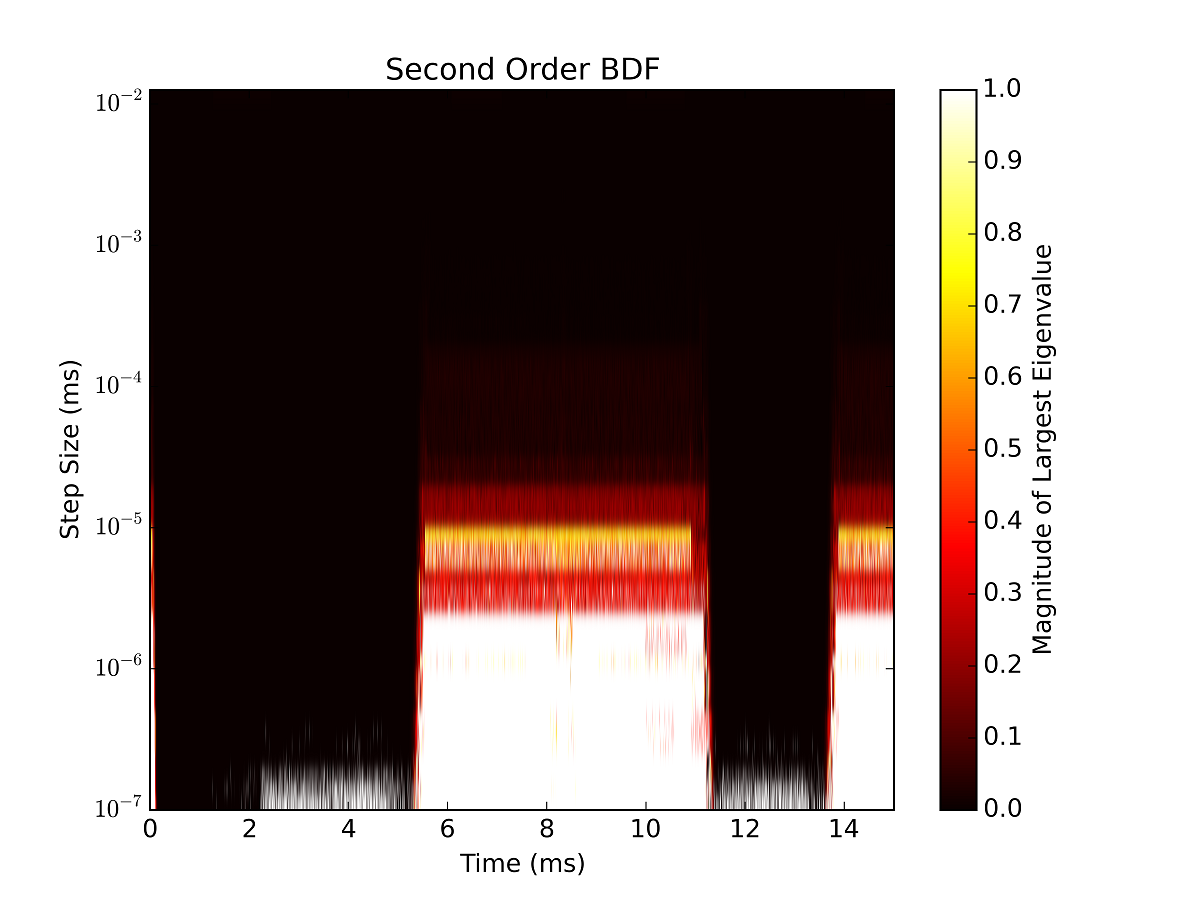}
	\caption{(left) The magnitude of the largest eigenvalue of Newton's method using a 1st order BDF with a step size indicated by the $y$-axis of the plot.
		(right) The magnitude of the largest eigenvalue using a 2nd order BDF.
		To advance in time, a default step size of $\Delta t = 10^{-3}$ ms was used in both cases.
		Note that white areas in the figure denote regions where the solver is believed to be unstable or where the solving process failed to converge after 20 iterations.
	}
	\label{fig:time step-sweep}	
\end{figure}

To highlight this impact and our method's ability to detect these issues prior to failure, we took the power channel example and integrated over a 15 ms window using a step size of $\Delta t = 10^{-3}$ ms.
At each of these reference steps, we then performed a sweep over 20 candidate step sizes ranging from $\Delta t = 10^{-2}$ ms down to $\Delta t = 10^{-7}$ ms for both a 1st and 2nd order BDF.
This sweep in $\Delta t$ has two key impacts on the nonlinear system that must be solved.
For both integrators,  $\Delta t$ shows up explicitly in the nonlinear system of equations.
However for the 2nd order BDF, the changes to $\Delta t$ implicitly appear in the coefficients of the solution history used by the method as the step size is no longer constant.

The results of this analysis are shown in Fig.~\ref{fig:time step-sweep}.
Despite the difference in the integration order, both integrators have a similar pattern: at certain intervals in time, decreasing the step size results in solver failure. 
Two regions in particular, when $t\in(5.5, 11.5)$ and $t > 14$, the solver appears to fail at much larger time steps than in other regions.
In both cases, however, there is a region where convergence anomalies exists (roughly from $\Delta t = 10^{-5}$ to $\Delta t=10^{-6}$) prior to solver failure.

In general, our approach monitors the efficiency of an implemented system-solver pairing relative to its idealized convergence behavior in order to improve simulation performance, confidence, and accuracy. 
In the first example, we held the solver fixed and would make changes to the system to correct the implementation errors we identified.
In this example, the system was held fixed, and the solver parameters were adjusted.
In particular, the analysis here reveals a trade off in step size selection: smaller step sizes will result in more accurate approximations of the underlying ODE/DAE, but after some point, can negatively impact the convergence rate of the associated solver.  

\section{Conclusions}
\label{sec:conclusions}

In this manuscript, we presented a pair of methods for identifying and diagnosing numerical issues in homotopy problems that negatively impact the solver performance.
Conceptually, both approaches appeal to dynamical systems theory, and are similar to the methods used to characterize the convergence rates of nonlinear solvers applied to an idealized system.
By using techniques from black-box bifurcation studies, we are able to apply the approach to practical systems, and identify issues that result from a mismatch between the implemented system and solver and the idealized ones.

We demonstrated our approaches on examples taken from the simulation of electrical circuits.
Our first set of examples were on the relatively simple diode bridge with one type of error at a time: either implementation errors in some of the electrical components or system-level issues arising from the removal of the minimum conductances.
The approach we proposed identified when these convergence issues arose, and in the cases with implementation errors, were able to isolate the components with the issue.
We also tested the method on a larger system meant to mimic a power channel on an aircraft that had multiple errors contained within different subsystems, and determined when the step size was becoming ``too small''.
In this more complex system, these techniques were able to separate the errors associated with each subsystem, and pinpoint the components where these errors were injected.

The current work has focused mainly on passive observation: detecting issues in the solver or underlying system and reporting these errors to the user.
However, one next step is adapting the solver based on these observations to correct these convergence anomalies and get closer to the theoretically expected level of performance.
Although these techniques cannot be applied to all systems (they are more intrusive than some commercial software will permit), they too can be thought of as the adaption of techniques like the recursive projection method~\cite{shroff1993stabilization} or Newton-Krylov-GMRES~\cite{kelley2003solving} to problems where the underlying system is another nonlinear solver.
Additionally, there are assumptions about the smoothness of the underlying system that are made implicitly throughout the manuscript.
The theory behind discontinuity induced bifurcations can explain some of the phenomena we have observed, such as sudden jumps in the eigenvalues, but at the current time, we have not implemented any checks targeted towards detecting the presence of discontinuities.
In principle this could provide useful information to the system and solver developers to improve system smoothness and solver selection.
In the end, we have demonstrated that even relatively simple techniques for the analysis of black-box dynamical systems can provide useful information about the practical performance of solvers when applied to actual systems.
Just as more sophisticated techniques for system analysis exist, the set of methods presented here are a starting point that could and should be augmented with more sophisticated approaches in the future.

\bibliographystyle{siam}
\bibliography{dynamicsrefs}

\end{document}